\documentclass[11pt]{article}
\catcode`\@=11 \catcode`\@=12

\usepackage{amsmath}
\usepackage{amssymb}
\usepackage{amsthm}
\usepackage{oldgerm}

\usepackage[pdftex,colorlinks,urlcolor=blue,pdfstartview=FitH]{hyperref}

\newcommand{\Rm}{\mathbb{R}}

\newcommand{\mF}{\ensuremath{\mathcal{F}}}

\newcommand{\Nm}{\ensuremath{\mathbb{N}}}
\newcommand{\Zm}{\ensuremath{\mathbb{Z}}}

\newcommand{\mA}{\ensuremath{\mathcal{A}}}

\newcommand{\mI}{\ensuremath{\mathcal{I}}}

\newcommand{\Tm}{\ensuremath{\mathbb{T}}}

\newtheorem{lem}{Lemma}
\newtheorem{ad}{Addendum}
\newtheorem{hyp}{Hypothesis}
\newtheorem{thm}{Theorem}

\newtheorem{cor}[lem]{Corollary}
\newtheorem{prop}[lem]{Proposition}
\newtheorem{property}[lem]{Property}
\newtheorem{defn}[lem]{Definition}

\def\proof {\noindent{\sc{Proof. }}}
\def\qed {\mbox{}\hfill {\small \fbox{}} \\}
\def\lto{\longrightarrow}
\def\lmto{\longmapsto}

\def\leq{\leqslant}
\def\geq{\geqslant}

\title{The Lax-Oleinik semi-group: a Hamiltonian point of view.}

\author{Patrick Bernard}

\date{CANPDE crash-course, Edimbourg, February 2011.
\footnote{Notes completed in May 2011, revised in january 2012.}}
\begin{document}
\maketitle

\begin{small}
The Weak KAM theory was developed by Fathi in order 
to study the dynamics of convex Hamiltonian systems.
It somehow makes a bridge between  viscosity solutions
of the Hamilton-Jacobi equation and Mather invariant sets
of Hamiltonian systems, although this was fully understood
only a posteriori.
These theories converge under 
the hypothesis of convexity, and the richness of applications 
mostly comes from this remarkable convergence.
In the present course, we provide an elementary exposition
of some of the basic concepts of weak KAM theory.
In a companion lecture, Albert Fathi exposes the aspects of his
theory which are more directly related to viscosity solutions.
Here on the contrary,  we focus on dynamical applications,
even if we also discuss some viscosity aspects
to underline the connections with Fathi's lecture.  
The fundamental reference on Weak KAM theory is the
still  unpublished book of Albert Fathi
\textit{Weak KAM theorem in Lagrangian dynamics}.
Although we do not offer new results, our exposition 
is original in several aspects. We only work with the Hamiltonian
and do not rely on the Lagrangian,
even if some proofs are directly
inspired from the classical Lagrangian proofs.
 This  approach is made easier
by the choice of a somewhat specific setting. We work on $\Rm^d$
and make  uniform hypotheses on the Hamiltonian.
This allows us to replace some  compactness arguments by explicit
estimates. For the most interesting dynamical
applications however, the compactness of the configuration space
remains a useful hypothesis and we retrieve it by
considering periodic (in space) Hamiltonians.
Our exposition is centered on the Cauchy problem for the 
Hamilton-Jacobi equation and the Lax-Oleinik evolution 
operators associated to it.
Dynamical applications are reached by considering fixed
points of these evolution operators, the Weak KAM solutions.
The evolution operators can also be used for their
regularizing properties, this opens a second way to dynamical
applications.

\tableofcontents

\end{small}

\section{The method of characteristics, existence and uniqueness of regular solutions.}\label{L1}
We consider a $C^2$ Hamiltonian
$$H(t,q,p):\Rm\times \Rm^d \times \Rm^{d*}\lto \Rm
$$
and study the associated Hamiltonian system
\begin{equation}\label{HS}\tag{HS}
\dot q(t)=\partial _pH(t,q(t),p(t))\quad,\quad
\dot p(t)=-\partial _qH(t,q(t),p(t))
\end{equation}
and Hamilton-Jacobi equation
\begin{equation}\label{HJ}\tag{HJ}
\partial_t u+H(t,q,\partial_qu(t,q))=0.
\end{equation}
We denote by $X_H(x)=X_H(q,p)$
the Hamiltonian vector field 
$X_H=JdH$, where $J$ is the matrix
$$
J=\left[\begin{array}{cc}
0&I\\-I&0
\end{array}\right].
$$
The Hamiltonian system can be written in condensed
terms $\dot x(t)=X_H(t,x(t))$.
We will always assume that that the solutions extend to $\Rm$.
We denote by
 $$
\varphi_{\tau}^t=(Q_{\tau}^t,P_{\tau}^t):\Rm^d \times \Rm^{d*}\lto \Rm^d \times \Rm^{d*}
$$
the flow map which, to a point $x\in {\mathbf{T}}^*\Rm^d$, associate
the value at time $t$ of the solution $x(s)$ of (\ref{HS})
which satisfies $x(\tau)=x$.

If $u(t,q)$ solves (\ref{HJ}), and if $q(s)$
is a curve in $\Rm^d$, then the formula
\begin{equation}\label{deltau}
u(t_1,q(t_1))-u(t_0,q(t_0))=
\int_{t_0}^{t_1}\partial_q u(s,q(s))\cdot
 \dot q(s)-
H(s, \partial_qu(s,q(s)))ds
 \end{equation}
follows from an obvious computation.
The integral on the right hand side is the
\textbf{Hamiltonian action} of the curve
$s\lmto (q(s),\partial_q u(s,q(s))).$
The  \textbf{Hamiltonian action} of the 
curve $(q(s),p(s))$ on the interval
$[t_0,t_1]$ is the quantity
$$
\int_{t_0}^{t_1}
p(s) \cdot  \dot q(s)-H(s,q(s),p(s))ds.
$$
A classical and important property of the Hamiltonian actions
is that orbits are critical points of this functional.
More precisely, we have:

\begin{prop}\label{critical}
The $C^2$ curve
$x(t)=(q(t),p(t)):[t_0,t_1]\lto \Rm^d \times \Rm^{d*}$
solves (\ref{HS}) if and only if the equality
$$
\left. \frac{d}{ds} \right|_{s=0}\left(
\int _{t_0}^{t_1}p(t,s)\cdot \dot q(t,s)-H(t,q(t,s),p(t,s))dt
\right) =0,
$$
where the dot is the derivative with respect to $t$,
holds for each $C^2$  variation 
$x(t,s)=(q(t,s),p(t,s)):[t_0,t_1]\times \Rm\lto \Rm^d \times \Rm^{d*}$ 
fixing the endpoints, which means  that
$x(t,0)=x(t)$ for each $t$ and that $q(t_0,s)=q(t_0)$ and 
$q(t_1,s)=q(t_1)$ for each $s$.
\end{prop}

\proof
We set $\theta(t)=\partial_sq(t,0)$,
 $\zeta(t)=\partial_sp(t,0)$ and compute:
\begin{align*}
&\left. \frac{d}{ds} \right|_{s=0}\left(
\int _{t_0}^{t_1}p(t,s)\dot q(t,s)-H(t,q(t,s),p(t,s))dt
\right) \\
&= \int _{t_0}^{t_1} p(t)\dot \theta(t)+
\zeta(t)\dot q(t)-\partial_q H(t,q(t),p(t) )\theta(t)
-\partial_p H(t,q(t),p(t))\zeta(t) dt\\
&= p(t_1)\theta(t_1)-p(t_0)\theta(t_0)+\int _{t_0}^{t_1}\big(
\dot q(t)-\partial_pH(t,q(t),p(t))
\big)\zeta(t)dt\\
&-\int _{t_0}^{t_1}
\big(
\dot p(t)+\partial_qH(t,q(t),p(t))
\big)\theta(t)dt.
\end{align*}
As a consequence, the derivative of the action vanishes
if $(q(t),p(t))$ is a Hamiltonian trajectory and if the variation $q(t,s)$
is fixing  the boundaries.
Conversely,  this computation can be applied 
to the variation $q(t,s)=q(t)+s\theta(t), p(t,s)=p(t)+s\zeta(t)$,
and implies that
$$\int _{t_0}^{t_1}\big(
\dot q(t)-\partial_pH(t,q(t),p(t))
\big)\zeta(t)dt-\int _{t_0}^{t_1}
\big(
\dot p(t)+\partial_qH(t,q(t),p(t))
\big)\theta(t)dt=0
$$
for each  $C^2$ curve $\theta(t)$ vanishing on the boundary
and each  $C^2$ curve $\zeta(t)$. This implies that 
$
\dot q(t)-\partial_pH(t,q(t),p(t))\equiv 0
$ and 
$
\dot p(t)+\partial_qH(t,q(t),p(t))\equiv 0
$.
\qed

We now return to 
the connections between
(\ref{HS}) and (\ref{HJ}). A function is said of class $C^{1,1}$
if it is differentiable and if its differential is Lipschitz. It is 
said of class $C^{1,1}_{loc}$ if it is differentiable with a 
locally Lipschitz differential.
The Theorem of Rademacher states that a locally Lipschitz function 
is differentiable almost everywhere. 

\begin{thm}\label{thmchar1}
Let $\Omega\subset\Rm\times \Rm^d$
be an open set, and let $u(t,q):\Omega\lto \Rm$
be a $C^{1,1}_{loc}$  solution of  (\ref{HJ}).
Let $q(t):[t_0,t_1]\lto \Rm^d$ be a $C^1$
curve such that $(t,q(t))\in \Omega$ 
and 
$$\dot q(t)=\partial _pH\big(t,q(t),\partial_q u(t,q(t))
\big)
$$
 for each
$t\in [t_0,t_1]$.
Then, setting $p(t)=\partial_q u(t,q(t))$,
the curve $(q(t),p(t))$ is $C^1$ and it solves (\ref{HS}).
\end{thm}
The curves $q(t)$ satisfying the hypothesis
of the theorem, as  well as the associated
trajectories $(q(t),p(t))$ are called the
\textbf{characteristics} of $u$.

\proof
Let $\theta(t):[t_0,t_1]\lto \Rm^d$ be a smooth curve
vanishing on the boundaries.
We define  $q(t,s):= q(t)+s\theta(t)$
and 
$
p(t,s):= \partial_qu(t,q(t,s)))
$.
Our hypothesis is that 
$\dot q(t)=\partial_p H(t,q(t),p(t))$,
which is the first part of (\ref{HS}).
For each $s$, we have
$$
u(t_1,q(t_1))-u(t_0,q(t_0))
=\int_{t_0}^{t_1}
p(t,s)\cdot\dot q(t,s)-H(t,q(t,s),p(t,s))ds
$$ 
hence $\left. \frac{d}{ds} \right|_{s=0}\left(
\int _{t_0}^{t_1}p(t,s)\dot q(t,s)-H(t,q(t,s),p(t,s))dt
\right)=0$.
We now claim that 
$$
\int _{t_0}^{t_1}
\partial_qH(t,q(t),p(t))
\cdot \theta(t) -p(t)\dot \theta(t) \, dt
=\left. \frac{d}{ds} \right|_{s=0}\left(
\int _{t_0}^{t_1}p(t,s)\dot q(t,s)-H(t,q(t,s),p(t,s))dt
\right).
$$
Assuming the claim, we obtain the equality 
$\int _{t_0}^{t_1}
\partial_qH(t,q(t),p(t))
\cdot \theta(t) -p(t)\cdot \dot \theta(t) \, dt
=0$ for each smooth function $\theta$ vanishing at the boundary.
In other words, we have 
$$\dot p(t)=-\partial_q H(t,q(t),p(t))
$$
in the sense of distributions.
Since the right hand side is continuous, this implies that $p$ is $C^1$ and that
the equality holds for each $t$. We have proved the theorem, assuming the claim.

The claim can be proved by an easy  computation in the case where $u$ is $C^2$.
Under the assumption  that $u$ is only $C^{1,1}_{loc}$, the map $p$ is only 
locally Lipschitz,
and some care is necessary.
For each fixed $\theta$,  we have 
\begin{align*}
\partial_qH(t,q(t,s),p(t,s))
\cdot \theta(t) -p(t,s)\cdot \dot \theta(t) 
&=
\partial_qH(t,q(t),p(t))
\cdot \theta(t) -p(t)\cdot \dot \theta(t)+O(s)\\
\partial_t q(t,s)-\partial_pH(t,q(t,s),p(t,s))
&=
\dot q -\partial_pH(t,q(t),p(t))+O(s)=O(s)
\end{align*}
where $O(s)$ is uniform in $t$. We then have, for small $S>0$,
\begin{align*}
&\int _{t_0}^{t_1}
\partial_qH(t,q(t),p(t))
\cdot \theta(t) -p(t)\cdot \dot \theta(t) \, dt\\
=&O(S)+\frac{1}{S}\int _{t_0}^{t_1}\int_0^S
\partial_qH(t,q(t,s),p(t,s))
\cdot \theta(t) -p(t,s)\cdot \dot \theta(t) \,ds dt \\
=&O(S)+\frac{1}{S}\int _{t_0}^{t_1}\int_0^S
\partial_qH
\cdot \partial_s q -p\cdot \partial_{st}q
+\big(\partial_t q-\partial_pH
\big)\cdot \partial_sp\;ds dt \\
=& O(S)+\frac{1}{S}\int _{t_0}^{t_1}\big[
p\cdot \partial_t q -H\big]_0^S\, dt
=
O(S)+\frac{1}{S}\left[\int _{t_0}^{t_1}
p\cdot \partial_t q -H\, dt\right]_0^S
.
\end{align*}
We obtain the claimed equality 
at the limit  $S\lto 0$. 
\qed

The following restatement of Theorem 
\ref{thmchar1} has a
more geometric flavor:

\begin{cor}\label{thmchar2}
Let $\Omega\subset\Rm\times \Rm^d$
be an open set, and let $u(t,q):\Omega\lto \Rm$
be a $C^{1,1}_{loc}$  solution of the Hamilton Jacobi equation (\ref{HJ}).
Then the extended Hamiltonian vector-field 
$Y_H=(1,X_H)$ is tangent to the graph
$$
G:=\{(t,q,\partial_q u): (t,q)\in \Omega\}.
$$
\end{cor}

\proof
Let us fix a point $(t_0,q_0)$ in $\Omega$.
By the Cauchy-Lipschitz theorem, there exists 
a solution $q(t)$ of the ordinary differential equation 
$\dot q =\partial_pH\big(t,q(t),\partial_qu(t,q(t))\big)$, 
defined on an open time interval containing $t_0$ and such that $q(t_0)=q_0$.
Let us define as above $p(t):=\partial_qu(t,q(t))$. 
The curve $(t,q(t),p(t))$ is contained in the graph $G$, and
we deduce from Theorem \ref{thmchar1} that it solves (\ref{HS}).
As a consequence, the derivative $Y_H$ of the curve $(t,q(t),p(t))$ 
is tangent to $G$.
\qed

\begin{cor}\label{flow}
Let $u(t,q)$ be a $C^{1,1}_{loc}$ solution of (\ref{HJ}) defined on the open set 
$\Omega=]t_0,t_1[\times \Rm^d$. Then, for each $s$ and $t$ in $]t_0,t_1[$
we have 
$$
\Gamma_t=\varphi_s^t(\Gamma_s),
$$
where $\Gamma_t$ is defined  by
$$
\Gamma_t:=\{(q,du_t(q)): q\in \Rm^d\}.
$$
\end{cor}
\proof
Let $(q_s,p_s)$ be a point in $\Gamma_s$.
Let us consider the Lipschitz map 
$$F(t,q):= \partial_pH(t,q,\partial_qu(t,q)),
$$
and consider the differential equation 
$
\dot q(t)=F(t,q(t)).
$
By the Cauchy-Peano Theorem, there exists a solution 
$q(t)$ of this equation, defined on the   interval 
$]t^-,t^+[\ni s$, and such that $q(s)=q_s$.
Setting $p(t)= \partial_qu(t,q(t))$, 
Theorem \ref{thmchar1}
implies that the curve $(q(t),p(t))$
solves (\ref{HS}).
We can chose $t^+$ such that either $t^+=t^1$
or the curve $q(t)$ is unbounded on $[s,t^+[$.
The second case is not possible because 
$(q(t),p(t))$ is a solution of (\ref{HS}),
which is complete, hence 
we can take $t^+=t_1$.
Similarly, we can take $t^-=t_0$.
We have proved that 
$(q(t),p(t))$ is the Hamiltonian orbit of the point
$(q_s,p_s)$.
Then, for each $t\in ]t_0,t_1[$, we have 
$$\varphi_s^t(q_s,p_s)=(q(t),p(t))=
(q(t),\partial_qu(t,q(t)))\in \Gamma_t.
$$
Since this holds for each $(q_s,p_s) \in \Gamma_s$,
we conclude that $\varphi_s^t(\Gamma_s)\subset \Gamma_t$
for each $s,t\in ]t_0,t_1[$.
By symmetry, this inclusion is an equality.
\qed

Let us now consider  an initial condition 
$u_{0}(q)$ and study the Cauchy problem consisting of finding
a  solution $u(t,q)$ of (\ref{HJ}) such that 
$u(0,q)=u_{0}(q)$.

\begin{prop}
Given a time interval $]t_0,t_1[$ containing the initial time $t=0$
and a $C^{1,1}_{loc}$ initial condition $u_{0}$, there is at most
one $C^{1,1}_{loc}$ solution $u(t,q):]t_0,t_1[\times \Rm^d$ of (\ref{HJ})
such that $u(0,q)=u_{0}(q)$ for all $q\in \Rm^d$.
\end{prop}

\proof
Let $u$ and $\tilde u$ be two solutions of this Cauchy problem.
Let us associate to them the graphs $\Gamma_t$ and $\tilde \Gamma_t$, $t\in ]t_0,t_1[$.
Since $\tilde u(\tau,q)=u(\tau,q)$, we have $\Gamma_{\tau}=\tilde \Gamma_{\tau}$ hence,
by Corollary \ref{flow}, 
$$
\Gamma_t=\varphi_{\tau}^t \big(\Gamma_{\tau}\big) =
\varphi_{\tau}^t \big(\tilde \Gamma_{\tau}\big) =
\tilde \Gamma_t.
$$
We conclude that $\partial_qu=\partial_q\tilde u$, and then, from (\ref{HJ}),
that $\partial_tu=\partial_t\tilde u$. 
The functions $u$ and $\tilde u$ thus have the same differential on $]t_0,t_1[\times \Rm^d$,
hence they differ by a constant.
Finally, since these functions have the same value on $\{\tau\}\times \Rm^d$, 
they are equal. 
\qed

To study the existence problem, 
 we lift the function $u_0$ to the surface $\Gamma_0$
by defining $w_0=u_0\circ \pi$, where $\pi$ is the projection 
$(q,p)\lmto q$ (later we will also use the symbol $\pi$  to denote the projection 
$(t,q,p)\lmto (t,q)$).
It is then useful to work in a more general setting:

A \textbf{geometric initial condition}
is the data of a subset $\Gamma_0\subset \Rm^d\times \Rm^{d*}$
and of a function $w_0:\Gamma_0\lto \Rm$ such that 
$dw_0=pdq$ on $\Gamma_0$. More precisely, we require that the equality 
$\partial_s(w_0(q(s),p(s)))=p(s)\partial_s q(s)$ holds almost everywhere for each 
Lipschitz curve $(q(s),p(s))$ on $\Gamma_0$. 
We will consider mainly two types of geometric initial conditions:

\begin{itemize}
\item The geometric initial condition $(\Gamma_0,w_0=u_0\circ \pi)$ associated 
to the $C^1$ initial condition $u_0$. 
\item The geometric initial condition $(\Gamma_0=\{q_0\}\times \Rm^{d*}, w_0=0)$, for $q_0\in \Rm^d$.
\end{itemize}

Given the geometric initial condition $(\Gamma_0,w_0)$, we define:
\begin{equation}\label{G}\tag{G}
G:=\bigcup_{t\in ]t_0,t_1[}\{t\}\times \varphi_{0}^t
\big( \Gamma_{0}\big)
\end{equation}
and, denoting by $\dot Q_t^s(x)$ the derivative with respect to $s$, the function 
\begin{align}\notag
w:\quad G&\lto \Rm \\\label{w} \tag{w}
(t,x)&\lmto 
w_0(Q_t^{0}(x))
+\int_{0}^t P_t^s(x) \dot Q_t^s(x) -H(s,\varphi_{t}^s(x))
ds.
\end{align}
The pair $(G,w)$ is called the \textbf{geometric solution} emanating from 
the geometric initial condition $(\Gamma_0,w_0)$.

This definition is motivated by the following observation:
Assume that a $C^2$ solution $u(t,q)$ of (\ref{HJ})   emanating from
the genuine initial condition $u_0$ exists.
Let  $(\Gamma_0,w_0)$ be  the geometric
initial condition associated to $u_0$.
Let $G$ be the graph of $\partial_qu$, as defined in  Corollary \ref{flow},
and let $w$ be the function defined on $G$ by 
$w:= u\circ \pi$.
Then, $(G,w)$ is the geometric solution emanating from the geometric initial condition 
$\Gamma_0$.
This follows immediately from Corollary \ref{flow} and equation (\ref{deltau}).
In general, we have:
\begin{prop}\label{solution}
Let $(\Gamma_0,w_0)$ be a geometric initial condition, and let $(G,w)$ be the geometric 
solution emanating from  $(\Gamma_0,w_0)$.
Then, the function $w$  satisfies
 $dw=pdq-Hdt$ on $G$. More precisely, for each Lipschitz curve 
$Y(s)=(T(s),\theta(s),\zeta(s))$ contained in $G$, then for a. e. $s$, 
$$\frac{d}{ds} \big(w(T(s),\theta(s),\zeta(s))\big)=
\zeta(s)\frac{d\theta}{ds}-H(Y(s))
\frac{dT}{ds}.
$$
\end{prop}

\proof
Let us first consider a $C^{2}$
curve  $Y(s)=(T(s),\theta(s),\zeta(s))$ on $G$.
We  set $q(t,s)=Q_{T(s)}^t(\theta(s),\zeta(s))$
and
$p(t,s)=P_{T(s)}^t(\theta(s),\zeta(s))$,
and finally $x(t,s)=(q(t,s),p(t,s))$.
We have 
$$
w(T(s),\theta(s),\zeta(s))=
w_0(q(0,s),p(0,s))+\int_{0}^{T(s)} 
p(t,s)\dot q(t,s)-H(t,x(t,s))dt.
$$
Since $dw_0=pdq$ on $\Gamma_0$,
the calculations in the proof of Proposition \ref{critical}
imply that 
\begin{align*}
\frac{d}{ds} (w\circ Y)&= p(0,s)\cdot\partial_sq(0,s)+
p(T(s),s)\cdot\partial_{s} q(T(s),s)
-p(0, s)\cdot\partial_ {s} q(0,s)\\
&+\Big(
p(T(s),s)\cdot \partial_t q(T(s),s)
-H(T(s),x(T(s),s))\Big)
\frac{dT}{ds}\\
&= \zeta(s)\left(
\partial_{s}q(T(s),s)+
\partial_t q(T(s),s)\frac{dT}{ds}
\right)+
H(Y(s))\frac{dT}{ds}.
\end{align*}
The desired equality follows from the observation that 
$d\theta/ds=\partial_t q(T(s),s)(dT/ds)+\partial_s q(T(s),s)$,
which can be seen by differentiating the equality $\theta(s)=q(T(s),s)$.

These computations, however, can't be applied directly  in the case where $Y(s)$
is only $C^1$, or, even worse, Lipschitz.
In this case, we will prove the desired equality in integral form
$$
[w\circ Y]_{S_0}^{S_1}=\int_{S_0}^{S_1} 
\zeta(s)\cdot \partial_s \theta (s)
-H\circ Y(s) \cdot\partial_s T(s) ds
$$
for each $S_0<S_1$. Fixing $S_0$ and $S_1$, we can approximate uniformly
the curve $Y(s)$ by a sequence 
$Y_n(s):[S_0,S_1]\lto \Rm\times \Rm^d\times\Rm^{d*}
$ 
of equi-Lipschitz smooth curves
such that $Y_n(S_0)=Y(S_0)$ and $Y_n(S_1)=Y(S_1)$.
To the curves $Y_n$, we associate $x_n(t,s)=(p_n(t,s),q_n(t,s))$ as above.
The functions  $x_n$ are equi-Lipschitz and converge uniformly to 
$x$.
In general, we don't have $Y_n(s)\in G$ on $]S_0,S_1[$, hence we don't
have $x_n(0,s)\in \Gamma_0$, and we cannot express $\partial_s w(x_n(0,s))$
as we did  above. Since this  is the only part of the above computation 
which used the inclusion $Y(s)\in G$, we can still get:
$$
\frac{d}{ds} (w\circ Y_n)= \frac{d}{ds} (w_0(x_n(0,s))
-p_n(0, s)\cdot\partial_ {s} q_n(0,s)+\zeta_n(s)\partial_s \theta_n(s)
+H(Y_n(s))\partial_sT_n(s).
$$
Noticing that 
$[w\circ Y]_{S_0}^{S_1}=[w\circ Y_n]_{S_0}^{S_1}$
and that 
$[w_0(x(0,.))]_{S_0}^{S_1}=[w_0(x_n(0,.))]_{S_0}^{S_1}$, we obtain 
\begin{align*}
[w\circ Y]_{S_0}^{S_1}
&=[w_0(x(0,.))]_{S_0}^{S_1}+
\int_{S_0}^{S_1}
-p_n(0, s)\cdot\partial_ {s} q_n(0,s)+\zeta_n(s)\partial_s \theta_n(s)
+H(Y_n(s))\partial_sT_n(s) ds\\
&=\int_{S_0}^{S_1} p(0,s)\cdot \partial_s q(0,s) -
p_n(0, s)\cdot\partial_{s} q_n(0,s) ds\\
&+\int_{S_0}^{S_1}\zeta_n(s)\partial_s \theta_n(s)
+H(Y_n(s))\partial_sT_n(s) ds.
\end{align*}
We derive the desired formula at the limit  $n\lto \infty$, along a subsequence
such that 
$$
\partial_sq_n(0,.)\rightharpoonup \partial_s q(0,.),\quad 
\partial_s\theta_n \rightharpoonup \partial_s \theta, \quad
\partial_sT_n\rightharpoonup\partial_sT
$$
weakly-$\star$ in $L^{\infty}$, taken into account  that 
$$
p_n(0,.)\lto p(0,.),\quad \zeta_n(s)\lto \zeta(s), \quad H(Y_n(s))\lto H(Y(s))
$$
uniformly, hence strongly in $L^1$. Recall  that a sequence of curves 
$f_n:[t_0,t_1]\lto \Rm^d$ is said to converge to $f$
 weakly-$\star$ in $L^{\infty}$
if $\int_{t_0}^{t_1}f_ng dt\lto\int_{t_0}^{t_1}fgdt$ 
 for each $L^1$ curve $g:[t_0,t_1]\lto \Rm^d$.
We have used two classical properties of the weak-$\star$
convergence:
\begin{itemize}
 \item A uniformly bounded sequence of functions has a subsequence which
 has a weak-$\star$
limit.
\item The convergence  $\int_{t_0}^{t_1}f_ng_n dt \lto \int fg dt$ 
holds if $f_n\rightharpoonup f$ weakly-$\star$ in $L^{\infty}$ and if 
$g_n\rightarrow g$  strongly in $L^1$.
\hfill \fbox{}
\end{itemize}

\begin{cor}\label{ex}
If there exists a locally Lipschitz map $\chi:\Omega \lto \Rm^{d*}$
on some open subset $\Omega$ of $]t_0,t_1[\times \Rm^d$
such that $(t,q,\chi(t,q))\subset G$ for all 
$(t,q)\in \Omega$, then the function 
$$u(t,q):=w(t,q, \chi(t,q))
$$
is $C^1$ and it solves (\ref{HJ}) on $\Omega$. Moreover, we have $\partial_qu=\chi$.
\end{cor}
\proof
For each $C^1$ curve $(T(s),Q(s))$ in $\Omega$,
the curve 
$$
Y(s)=(T(s),Q(s),\chi(T(s),Q(s))
$$
 is Lipschitz, hence,
by Proposition \ref{solution},
we have 
\begin{align*}
\partial_s u(T(s),Q(s))
&=\partial_sw\big(T(s),Q(s),\chi(T(s),Q(s))\big) \\
&=\chi\big(T(s),Q(s)\big) \cdot \partial_s Q(s)-
H\big(T(s),Q(s),\chi(T(s),Q(s)\big) \partial_sT(s)
\end{align*}
almost everywhere.
Since the right hand side in this expression is continuous,
we conclude that the Lipschitz functions $u(T(s),Q(s))$ is actually
differentiable at each point, the equality above being satisfied 
everywhere. Since this holds for each $C^1$ curve $(T(s), Q(s))$, 
the function $u$ has to be  differentiable, with
$
\partial_qu(t,q)=\chi(t,q)$ and 
$
\partial_tu(t,q)+H(t,q,\chi(t,q))=0. 
$
\qed

We have reduced the existence problem to the study of the geometric solution $G$.
We need an additional hypothesis to obtain a local existence result.
We will rest on the following one,  which it is stronger than would really be necessary,
but will allow us to rest on simple estimates in this course.

\begin{hyp}\label{h1}
There exists a  constant $M$
such that 
$$
\|d^2H(t,q,p)\|\leq M
$$
for each $(t,q,p)$.
\end{hyp}

This hypothesis implies that the Hamiltonian
vector-field is Lipschitz, hence that the Hamiltonian flow
is complete.
The hypothesis can be exploited further to estimate the 
differential  
$$d\varphi_0^t=\left[
\begin{array}{cc}
\partial_q Q_0^t(x)&
\partial_p Q_0^t(x)\\
\partial_q P_0^t(x)&
\partial_p P_0^t(x)
\end{array}
\right]
$$
using the variational equation
$$
\left[
\begin{array}{cc}
\partial_q \dot Q_0^t(x)&
\partial_p \dot Q_0^t(x)\\
\partial_q \dot P_0^t(x)&
\partial_p \dot P_0^t(x)
\end{array}
\right]
=
\left[
\begin{array}{cc}
\partial_{qp}H(t,x)&
\partial_{pp}H(t,x)\\
-\partial_{qp}H(t,x)&
-\partial_{pp}H(t,x)
\end{array}
\right]
\left[
\begin{array}{cc}
\partial_q  Q_0^t(x)&
\partial_p Q_0^t(x)\\
\partial_q  P_0^t(x)&
\partial_p  P_0^t(x)
\end{array}
\right].
$$
We obtain the following estimates:
$$
\|d\varphi_{\tau}^t -I\|\leq e^{M|t-\tau|}-1
$$
which implies, for $|t-\tau|\leq 1/M$, that
\begin{equation}\label{M}\tag{M}
\|d\varphi_{\tau}^t -I\|\leq 2M|t-\tau|
\end{equation}
or componentwise (taking $\tau=0$, and assuming that $|t|\leq M$):
$$
\|\partial_q  Q_0^t-I\|\leq 2M|t|\quad,\quad
\|\partial_p P_0^t-I\|\leq 2M|t|\quad,\quad
\|\partial_q P_0^t\|\leq 2M|t|\quad,\quad
\|\partial_p Q_0^t\|\leq 2M|t|.
$$
We can now prove:

\begin{thm}\label{existence}
Let $H:\Rm\times \Rm^d\times (\Rm^d)^*$ be a $C^2$ Hamiltonian
satisfying 
 Hypothesis \ref{h1}.
Let $u_0$ be a  $C^{1,1}$ initial condition.
There exists
a time $T>0$ and a $C^{1,1}_{loc}$ solution 
$u(t,q):]-T,T[\times \Rm^d\lto \Rm$ of (\ref{HJ})
such that $u(0,q)=u_0(q)$.
Moreover, we can take 
$$
T=\Big(4M\big(1+Lip(du_0)\big)\Big)^{-1},
$$
and we have
$$
Lip(du_t)\leq Lip(du_0)+4|t|M\big(1+Lip(du_0)\big)^2
$$
when $|t|\leq T$. If the initial condition $u_0$ is $C^2$, then so is the solution
$u(t,q)$.
\end {thm}

\proof
Let $(\Gamma_0,w_0)$ be the geometric initial condition associated to $u_0$, and 
let $(G,w)$ be the geometric solution emanating from $(\Gamma_0,w_0)$.
We first prove that the restriction of $G$ to $]-T,T[\times \Rm^d$
is a graph.
It is enough to prove that the map
$$
 F(t,q):=\big(t, Q_{0}^t (q,du_0(q))\big)
$$
is a bi-Lipschitz homeomorphism of $]-T,T[\times \Rm^d$.
By  (\ref{M}), we have
$$
Lip(F-Id)\leq 2|t|M\big(1+Lip(du_0)\big)<1,
$$
provided 
$|t|< \Big(2M\big(1+Lip(du_0)\big)\Big)^{-1}.
$
We conclude using the classical Proposition \ref{global}
of the Appendix that $F$ is a bi-Lipschitz homeomorphism 
of $]-T,T[\times \Rm^d$. Moreover,
if $u_0$ is $C^2$, then $F$ is a $C^1$ diffeomorphism.
Since $F$ is a homeomorphism preserving $t$, we can denote by
 by $(t,Z(t,q))$ its inverse. By Proposition \ref{global}, we have
$$
Lip(Z)\leq \frac{1}{1-2|t|M(1+Lip(du_0))},
$$
and, under the assumption that $|t|\leq T$ (as defined in the statement),
we obtain 
$$
Lip(Z)\leq 1+4M|t|(1+Lip(du_0)) \leq 2.
$$
We have just used here that $(1-a)^{-1}\leq 1+2a$ for $a\in[0,1/2]$.
We set
$$\chi(t,q)=P_0^t\big(Z(t,q),du_0(Z(t,q))\big),
$$
in such a way that  $G$ 
is the graph of $\chi$ on $]-T,T[\times \Rm^d$.
Observing that  $\chi$ is Lipschitz, 
 we conclude from Corollary \ref{ex} that the function  $u(t,q):=w(t,q,\chi(t,q))$
solves (\ref{HJ}). Moreover, we have $u(0,q)=u_0(q).$
 Corollary \ref{ex} also implies  that $du_t=\chi_t$ hence, in view of
 (\ref{M}), we have 
\begin{align*}
Lip(du_t)
&= Lip(\chi_t)\leq 2M|t| Lip(Z_t)+(1+2M|t|)Lip(du_0)Lip(Z_t)\\
&\leq 4M|t|+Lip(du_0)+Lip(du_0)(4M|t|(1+Lip(du_0)))+4M|t|Lip(du_0)\\
&\leq Lip(du_0)+4M|t|\big(1+Lip(du_0)\big)
\big(1+Lip(du_0)\big).
\end{align*}
\qed

\subsection{Exercise :}
Take $d=1$, 
$H(t,q,p)=(1/2)p^2$, and 
$u_0(q)=-q^2$, and prove that the $C^2$ solution can't be extended 
beyond $t=1/2$.

\section{Convexity, the twist property, and the generating function.}\label{sectionconvex}
We  make an  additional assumption on $H$. Once again, we make the assumption
in a stronger form than  would be necessary, this allows  to obtain simpler
statements:

\begin{hyp}\label{h2}
There exists $m>0$
such that 
$$
\partial^2_{pp} H\geq mId
$$
for each $(t,q,p)$, in the sense of quadratic forms.
\end{hyp}

Let us first study the consequences of this
hypothesis on the  structure of the flow.

\begin{prop}\label{twist}
There exists $\sigma>0$ such that the map 
$p\lmto Q_{0}^t(q,p)$ is $(mt/2)$-monotone
when $t\in ]0,\sigma]$, in the sense that the inequality
$$
(Q_{0}^t(q,p')-Q_{0}^t(q,p))\cdot (p'-p)
\geq mt|p'-p|^2/2
$$
holds for each $q\in \Rm^d$, each $t\in [0,\sigma]$.
As a consequence, it is a $C^1$ diffeomorphism onto $\Rm^d$.
\end{prop}
We say that the flow has the \textbf{Twist property}.

\proof
Fix  a point $q$ and denote by $F^t$
the map $p\lmto Q_0^t(q,p)$.
We have
$
dF^t(p)=\partial _pQ_0^t(q,p).
$
In order to estimate this linear map,
we recall the variational equation 
$$
\partial_p \dot Q_0^t(x)=
\partial_{qp}H(t,\varphi_0^t(x))\partial_p  Q_0^t(x)
+
\partial_{pp}H(t,\varphi_0^t(x))\partial_p  P_0^t(x).
$$
We deduce that 
$$
\partial_p \dot Q_0^t(x)-\partial^2_{pp}H(t,\varphi_0^t(x))=
\partial^2_{qp}H(t,\varphi_0^t(x))\partial_p  Q_0^t(x)
+
\partial^2_{pp}H(t,\varphi_0^t(x))(\partial_p  P_0^t(x)-Id)
$$
and then that 
$$
\|
\partial_p \dot Q_0^t(x)-\partial^2_{pp}H(t,\varphi_0^t(x))
\|\leq 2M^2t
$$
As a consequence, for $t\leq \sigma=m/(4M^2)$, 
we have
$$
\partial_p\dot Q_0^t\geq (m-2M^2t)I\geq (m/2) I
$$
in the sense of quadratic forms (note that the matrix 
$\partial_p\dot Q_0^t$ is not necessarily symmetric).
Since 
$$
\partial_pQ_0^t(x)=\int_0^t \partial_p\dot Q_0^s(x) ds,
$$
we conclude that 
$$
dF^t(p)=\partial_pQ_0^t(q,p)\geq (m/2)Id,
$$
which means that $(dF^t(p)z,z)\geq (m/2)|z|^2$
for each $z\in \Rm^{d*}$.
This estimate can be integrated, and implies the monotony of the map 
$F^t$:
\begin{align*}
(Q^t(q,p')-Q^t(q,p))\cdot (p'-p)
&=
\left(\int_0^1 \partial_pQ^t(q,p+s(p'-p))\cdot (p'-p)ds
\right)\cdot (p'-p)\\
&=
\int_0^1\left( \partial_pQ^t(q,p+s(p'-p))\cdot (p'-p)\right) ds\\
&\geq\int_0^1 (m/2)t(p'-p)\cdot (p'-p) ds
\geq (m/2)t(p'-p)\cdot (p'-p).
\end{align*}
It is then a classical result that the map $F^t$ is a $C^1$ diffeomorphism, see Proposition \ref{monotonemap} in the appendix.
\qed

\begin{cor}
The map 
$(t,q,p)\lmto (t,q, Q_0^t(q,p))$
is a $C^1$ diffeomorphism from 
$]0,\sigma[\times \Rm^d\times \Rm^{d*}$
onto its image
$]0,\sigma[\times \Rm^d\times \Rm^{d}$.
\end{cor}

We denote by $\rho_0(t,q_0,q_1)$ the unique momentum $p$
such that  $Q_0^t(t,q_0,\rho_0(t,q_0,q_1))=q_1$.
In other words, $\rho_0(t,q_0,q_1)$ is the initial 
momentum $p(0)$ of the unique orbit
$(q(s),p(s)):[0,t]\lto \Rm^d\times \Rm^{d*}$ of (\ref{HS})
which satisfies $q(0)=q_0$ and $q(t)=q_1$. 
By the Corollary we just proved, 
the map $\rho_0$ is $C^1$.
Similarly, we denote by $\rho_1(t,q_0,q_1)$ the unique momentum $p$
such that  $Q_t^0(t,q_1,\rho_1(t,q_0,q_1))=q_0$.
We can equivalently define   $\rho_1$ as 
$$
\rho_1(t,q_0,q_1)=P_0^t(t,q_0,\rho_0(t,q_0,q_1)).
$$
Considering the geometric initial condition 
$(\Gamma_0=\{q_0\}\times \Rm^{d*},w_0=0)$,
and the associated geometric solution $(G,w)$, we see that 
$$
G=\{(t,q,\rho_1(t,q_0,q)),\quad  (t,q)\in ]0,\sigma[\times \Rm^d\}.
$$
We conclude from Corollary \ref{ex} that there exists a genuine solution
of (\ref{HJ}) emanating from the geometric initial condition 
$(\{q_0\}\times \Rm^{d*},0)$.
We denote by $S^t(q_0,q)$ this solution. We have 
$$
S^t(q_0,q)=w(t,p, \rho_1(t,q_0,q))
$$
and 
$$\partial_q S^t(q_0,q)=\rho_1(t,q_0,q).
$$
In view of the definition of geometric solutions, 
the function $S$ can be written more explicitly
$$
S^t(q_0,q_1)=
\int_0^t
P_0^s(q_0,\rho_0(t,q_0,q_1))\dot Q_0^s(q_0,\rho_0(t,q_0,q_1))
-H(s, \varphi_0^s(q_0,\rho_0(t,q_0,q_1)) ds.
$$
In words, 
$S^t(q_0,q_1)$
is the action of the unique trajectory
 $(q(s),p(s)):[0,t]\lto \Rm^d\times \Rm^{d*}$ of (\ref{HS})
which satisfies $q(0)=q_0$ and $q(t)=q_1$.

We have defined the function $S^t(q_0,q_1)$
as the action of the unique orbit joining $q_0$
and $q_1$ between time $0$ and $t$. 
We can define similarly the function $S^t_{\tau}(q_0,q_1)$
as the action of the unique orbit joining 
$q_0$ to $q_1$ between time $\tau$ and time $t$, all this being well-defined provided $0<t-\tau<\sigma$.
It is possible to prove as above that 
the function $(s,q)\lmto S_s^{t}(q,q_1)$ solves the 
Hamilton-Jacobi equation 
$$
\partial_su +H(t,q,-\partial_qu)=0,
$$
on $s<t$, and that 
$$
\partial_qS^t(q,q_1)=\partial_qS_0^t(q,q_1)=-\rho_0(t,q,q_1).
$$
\textbf{Convention: }We shall from now on denote by $\partial_0S^t$ the partial differential with respect
to the first variable (which in our notations is often $q_0$),
 and by $\partial_1S^t$ the partial differential with respect
to the second  variable (which in our notations is often $q_1$).

The relations $\partial_0S=-\rho_0$, $\partial_1S=\rho_1$,
$\partial_tS=-H(t,q_1,\rho_1)=-H(0,q_0,\rho_0)$
that we have proved imply that the function $S$ is 
$C^2$. 
Moreover, since $\varphi_0^t(q_0,\rho_0(t,q_0,q_1))=(q_1,\rho_1(t,q_0,q_1)$, we have
$$
\varphi_0^t(q_0,-\partial_0S(q_0,q_1))=(q_1,\partial_1S^t(q_0,q_1)).
$$
We say that $S^t$ is a \textbf{generating function}
of the flow map $\varphi_0^t$.
See \cite{MDS}, chapter 9, for more material on generating functions.
 It is useful to estimate the  second differentials of $S$:

\begin{lem}\label{estimeesS}
The function $S$ is $C^2$ on $]0,\sigma[\times \Rm^d\times \Rm^d$, and the estimates
\begin{align*}
\partial^2_{00}S^t\geq \frac{c}{t} I\quad ,&\quad
\partial^2_{11}S^t \geq  \frac{c}{t} I\\
\|\partial^2_{00}S^t \| +
\|\partial^2_{01}S^t \|& +
 \|\partial^2_{01}S^t \|\leq \frac{C}{t}
\end{align*}
hold, with constants $c$ and $C$ which depend
only on $m$ and $M$.
\end{lem}

\proof
Let us first observe that 
$$
\partial^2_{11}S^t(q_0,q_1)=
\big(\partial_pP_0^t(q_0,\rho_0(t,q_0,q_1)\big)
\big(\partial_pQ_0^t(q_0,\rho_0(t,q_0,q_1))\big)^{-1},
$$
and  recall the estimates:
$$
\|\partial_pP_0^t-Id\|\leq 2Mt,\quad
\|\partial_pQ_0^t \|\leq 2Mt,\quad \partial_pQ_0^t\geq (mt/2)Id.
$$
We conclude  that (see Lemma \ref{a/b})
$$
(\partial_pQ_0^t)^{-1}\geq \frac{m}{8M^2t}Id
\quad, \quad 
\|(\partial_pQ_0^t)^{-1}\|\leq 2/(mt).
$$
Finally, we obtain that 
$$
\partial^2_{11}S(q_0,q_1)\geq
 \left( \frac{m}{8M^2t}-
\frac{4M}{m}\right) Id
\geq \frac{m}{16M^2t}Id
$$
provided $t\leq m^2/(64M^3)$.
The other estimates can be proved similarly, using the expressions
\begin{align*}
\partial^2_{00}S^t(q_0,q_1)&=-\big(\partial_pP_t^0(q_1,\rho_1(t,q_0,q_1)\big)
\big(\partial_pQ_t^0(q_1,\rho_1(t,q_0,q_1))\big)^{-1},\\
\partial^2_{10}S^t(q_0,q_1)&=\left(\partial_pQ_0^t
(q_0,\rho_0(t,q_0,p_0))\right)^{-1}.
\end{align*}
\qed

\begin{prop}\label{triangle}
Given times $t_1$ and $t_2$ such that $0<t_1<t_2<\sigma$, we have 
the triangle inequality
$$
S_0^{t_2}(q_0,q_2) 
\leq S_0^{t_1}(q_0,q_1)+S_{t_1}^{t_2}(q_1,q_2)
$$
for each $q_0, q_1,q_2$. 
Moreover, 
$S_0^{t_2}(q_0,q_2)=\min_q
\big(S_0^{t_1}(q_0,q)+S_{t_1}^{t_2}(q,q_2)\big).
$
\end{prop}

\proof
Let us consider the map
$$
q\lmto f(q)=S_0^{t_1}(q_0,q)+S_{t_1}^{t_2}(q,q_2).
$$
We have $d^2f\geq 2c$ hence the map $f$ is convex.
Now let us denote by $(q(s),p(s)):[0,t_2]\lto \Rm^d\times \Rm^{d*}$
the unique orbit which satisfies $q(0)=q_0$ and 
$q(t_2)=q_2$.
We can compute 
$$
df(q(t_1))=\partial_1S_0^{t_1}(q_0,q(t_1))+
\partial_0S_{t_1}^{t_2}(q(t_1),q_{2})
=p(t_1)-p(t_1)=0.
$$
The point $q(t_1)$ is thus a critical point of the convex function $f$, hence it is a minimum of this function.
We conclude that 
$$
S_0^{t_1}(q_0,q)+S_{t_1}^{t_2}(q,q_{2})
\geq 
S_0^{t_1}(q_0,q(t_1))+S_{t_1}^{t_2}(q(t_1),q_{2})
=S_0^{t_2}(q_0,q_{2})
$$
for all  $q$.
\qed
Under the convexity hypothesis  \ref{h2}, Theorem \ref{thmchar1}
can be extended to $C^1$ solutions:
\begin{thm}\label{thmchar3}
Let $\Omega\subset\Rm\times \Rm^d$
be an open set, and let $u(t,q):\Omega\lto \Rm$
be a $C^1$  solution of the Hamilton Jacobi equation (\ref{HJ}).
Let $q(t):[t_0,t_1]\lto \Rm^d$ be a $C^1$
curve such that $(t,q(t))\in \Omega$
and 
$$\dot q(t)=\partial _pH\big(q(t),\partial_q u(t,q(t))
\big)
$$
 for each
$t\in [t_0,t_1]$.
Then, setting $p(t)=\partial_q u(t,q(t))$,
the curve $(q(t),p(t))$ solves (\ref{HS}).
\end{thm}

\proof
As in the proof of Theorem \ref{thmchar1}, we consider 
 a  variation
$q(t,s)=q(t)+s\theta(t)$ of $q(t)$,
where $\theta$ is smooth and vanishes on the endpoints.
We choose the vertical variation 
$p(t,s)$ in such a way that the 
equation 
$$
\dot q(t,s)=\partial_pH(t,q(t,s),p(t,s))
$$
holds. 
The map $p(t,s)$ defined by this relation is differentiable in $s$, because $q$ and $\dot q$ are and because 
the matrix $\partial^2_{pp}H$ is invertible.
It is also useful to consider the other vertical variation 
$$
P(t,s):= \partial_qu(t,q(t,s)).
$$
Our hypothesis is that 
$\dot q(t)=\partial_p H(t,q(t),p(t))$,
which is the first part of (\ref{HS}).
We start as in the proof of Theorem  \ref{thmchar1}
with the following equality:
\begin{equation}\label{null}
\left. \frac{d}{ds} \right|_{s=0}\left(
\int _{t_0}^{t_1}p(t,s) \cdot \dot q(t,s)-H(t,q(t,s),p(t,s))dt
\right)=0.
\end{equation}
We deduce this equality from the observation that 
$s=0$ is a local minimum of the function 
$$
s\lmto F(s):=
\int _{t_0}^{t_1}p(t,s)\cdot \dot q(t,s)-H(t,q(t,s),p(t,s))dt
.
$$
This claim follows from the equality
$$
F(0)=u(t_1,q(t_1))-u(t_0,q(t_0))
=\int_{t_0}^{t_1}
P(t,s)\cdot\dot q(t,s)-H(t,q(t,s),P(t,s))ds,
$$ 
which holds for all  $s$, and from the inequality
$$
F(s)\geq \int_{t_0}^{t_1}
P(t,s)\cdot\dot q(t,s)-H(t,q(t,s),P(t,s))ds
$$
which results, in view of the convexity of $H$, from the computation
\begin{align*}
H(t,q(t,s),P(t,s))&
\geq (P(t,s)-p(t,s))\cdot \partial_pH(t,q(t,s),p(t,s)) +H(t,q(t,s),p(t,s))\\
& \geq (P(t,s)-p(t,s))\cdot \dot q(t,s)+H(t,q(t,s),p(t,s)).
\end{align*}
We have proved (\ref{null}). As in the proof of Theorem \ref{thmchar1},
we develop the left hand side and, after a simplification, we get
$$
\int_{t_0}^{t_1} p(t)\cdot \dot \theta(t)
-\partial_qH(t,q(t),p(t))\cdot \theta(t) \,dt =0.
$$
In other words, we have proved that $\dot p(t)=\partial_qH(t,q(t),p(t))$
in the sense of distributions. Since the right hand side is continuous, 
$p$ is $C^1$ and 
the equality holds in the genuine sense.
\qed
As in the $C^2$ case, we have the following corollary,
 see \cite{Fa:03}:

\begin{cor}\label{geomc1}
Let $u(t,q):]t_0,t_1[\times \Rm^d\lto \Rm$ be a $C^1$ solution of (\ref{HJ}). Then, for each $s$ and $t$ in $]t_0,t_1[$
we have 
$$
\Gamma_t=\varphi_s^t(\Gamma_s),
$$
where $\Gamma_t$ is defined  by
$$
\Gamma_t:=\{(q,du_t(q)): q\in \Rm^d\}.
$$
\end{cor}

\proof
This corollary follows from Theorem
\ref{thmchar3} in the same way as
 Corollary \ref{flow} follows from Theorem \ref{thmchar1}.
The only difference here is that the map
$$F(t,q):= \partial_pH(t,q,\partial_qu(t,q))
$$
is only continuous. By the Cauchy-Peano Theorem,
this is sufficient to imply the existence of solutions
to the associated differential equation, which is 
what we need to develop the argument.
\qed

A last property of the functions 
$S$ will be useful.
Assume that we are considering a family $H_{\mu},\mu \in I$
of Hamiltonians, where $I\subset \Rm$ is an interval,
 such that the whole function $H(\mu,t,q,p)$
is $C^2$ and such that each of the Hamiltonians
$H_{\mu}$ satisfy our hypotheses \ref{h1} and \ref{h2},
with uniform constants $m$ and $M$.
Then, for each value of $\mu$, we have the function 
$S^t(\mu;q_0,q_1)$, which is defined for 
$t\in ]0,\sigma]$, the bound $\sigma>0$ being independent of 
$\mu$.
Since everything we have done so far was based on the 
local inversion theorem, the function $S^t(\mu;q_0,q_1)$
is $C^1$ in $\mu$, or more precisely 
the function $(\mu,t,q_0,q_1)\lmto S^t(\mu;q_0,q_1)$
is $C^1$. 
Moreover, a computation similar to the proof
of Proposition \ref{critical} yields  
$$
\partial_{\mu}S^t(\mu;q_0,q_1)=-\int_0^t\partial_{\mu}H_{\mu}(s,q(\mu,s),p(\mu,s))ds,
$$
where $s\lmto (q(\mu,s),p(\mu,s))$ is the only $H_{\mu}$-trajectory 
satisfying $q(\mu,0)=q_0$ and $q(\mu,t)=q_1$.
We can exploit this remark when $H_{\mu}$
is the linear interpolation
$H_{\mu}=H_0+\mu(H_1-H_0)$ between two Hamiltonians
$H_0$ and $H_1$, 
and conclude the important monotony property:
\begin{equation}\tag{Monotone}\label{Monotone}
H_0\leq H_1 \quad \Rightarrow \quad 
S^t(0;q,q')\geq S^t(1;q,q').
\end{equation}

\subsection{Exercise :}
If $H(t,q,p)=h(p)$ is a function of $p$, then 
$$
S^t(q_0,q_1)=th^*\left( \frac{q_1-q_0}{t}\right),
$$
where $h^*$ is the Legendre transform of $h$.
As an example, when $H(t,q,p)=a|p|^2/2$, then
$$
S^t(q_0,q_1)=\frac{1}{2ta}|q_1-q_0|^2.
$$

\section{Extension of the generating function: The minimal action.}

A classical problem consists in finding an orbit 
$(q(t),p(t))$ of the Hamiltonian system such that 
$q(t_0)=q_0$ and $q(t_1)=q_1$, for given
$[t_0,t_1]\subset\Rm$, $q_0,q_1\in \Rm^d$. 
We have seen, under Hypotheses \ref{h1} and \ref{h2},
that this problem has a unique solution 
provided $t_0<t_1<t_0+\sigma$, where $\sigma$ is 
a constant depending only on $m$ an $M$.
The situation is more subtle for larger values of 
$t_1-t_0$.
In order to study it, it is useful to consider 
the function
$$
\mathfrak{S}:(\theta_1,\ldots,\theta_{n-1})
\lmto
S^{t/n}_0(q_0,\theta_1)+S^{2t/n}_{t/n}(\theta_1,\theta_2)+\cdots+
S^{t}_{(n-1)t/n}(\theta_{n-1},q_1),
$$
where we have taken $t_0=0$ and $t_1=t$ to simplify
notations, and where $n$ is an integer such that 
$t/n\leq \sigma$.
The critical points of $\mathfrak{S}$ are in one to one 
correspondence with the solutions of our problem:

\begin{lem}
The point $(\theta_1,\ldots, \theta_{n-1})$ is a critical point of $\mathfrak{S}$ if and only if there exists an 
orbit 
$(q(s),p(s)):[0,t]\lto \Rm^d\times \Rm^{d*}$
such that $q(0)=q_0$, $q(t)=q_1$, and 
$q(it/n)=\theta_i$ for $i=1,\ldots, n-1$.
This orbit is then unique, and its action is
$\mathfrak{S}(\theta_1,\ldots, \theta_{n-1})$.
\end{lem}

\proof
Let $(q(s),p(s))$ be the piecewise orbit defined
on $[it/n,(i+1)t/n]$ by the constraints $q(it/n)=\theta_i$
and $q((i+1)t/n)=\theta_{i+1}$.
The action of this piecewise orbit is 
$\mathfrak{S}(\theta_1,\cdots, \theta_{n-1})$.
The statement follows from the simple computation
$$
\partial_{\theta_i}
\mathfrak{S}=\partial_1S^{t/n}(\theta_{i-1},\theta_i)+
\partial_0S^{t/n}(\theta_i,\theta_{i+1})
=p^-(it/n)-p^+(it/n).
$$
\qed
Using this finite dimensional variational functional
is usually called the method of broken geodesics, see 
\cite{Ch}.
The function $\mathfrak{S}$ can be minimized
under  additional assumptions, for example:
\begin{hyp}\label{h3}
$$
\frac{m}{2}|p|^2-M\leq H(t,q,p)\leq \frac{M}{2}|p|^2+M.
$$
\end{hyp}
By exploiting the monotony property
(\ref{Monotone}), this hypothesis implies  that 
$$
\frac{1}{2tM}|q_1-q_0|^2-Mt
\leq S^t(q_0,q_1)\leq \frac{1}{2tm}|q_1-q_0|^2+Mt,
$$
and then that 
$$
\mathfrak{S}(\theta_1,\ldots,\theta_{n-1})
\geq \frac{n}{2tM}(|\theta_1-q_0|^2+
|\theta_2-\theta_1|^2+\cdots +|q_1-\theta_{n-1}|^2)
-Mt.
$$
As a consequence, the function $\mathfrak{S}$
is coercive and $C^2$, hence it has a minimum.
Notice that, although $\mathfrak{S}$ is convex separately
in each of its variables, it is not jointly convex.
It can have critical points which are not minima,
and it can have several different minima.
We denote by $A^t$ the value function
\begin{equation}\label{A}\tag{A}
A^t(q_0,q_1)=\min \mathfrak{S}=\min _{\theta_1,\theta_2,\ldots \theta_{n-1}}
\big(S^{t/n}_0(q_0,\theta_1)+S^{2t/n}_{t/n}(\theta_1,\theta_2)+
S^{t}_{(n-1)t/n}(\theta_{n-1},q_1)\big)
\end{equation}
where $n$ is any integer such that $t/n<\sigma$. 
The functions $A_{\tau}^t(q_0,q_1)$ are defined similarly for each $t\geq \tau.$
This notation is legitimate in view of the following:

\begin{lem}\label{estimeA}
The value of $A^t$ does 
 not depend on $n$ provided $t/n<\sigma$.
Moreover, we have
$$\frac{1}{2Mt}|q_1-q_0|^2-Mt\leq A^t(q_0,q_1)
\leq \frac{1}{2tm}|q_1-q_0|^2+Mt.
$$
\end{lem}
This statement implies that $A^t=S^t$ when $t<\sigma$ :
$A^t$ can be seen as an extension of $S^t$ beyond 
$t=\sigma$.
 
\proof
Since we have not yet proved the independence of $n$,
we temporarily denote by $A^t(q_0,q_1;n)$ the value of the minimum.
We have 
\begin{align*}
A^t(q_0,q_1;n)
&\geq \min_{\theta_1,\theta_2,\ldots \theta_{n-1}}
\left(\frac{n}{2Mt}\big(|\theta_1-q_0|^2
+\cdots+|q_1-\theta_{n-1}|^2\big)-Mt\right)\\
&=\frac{1}{2Mt}|q_1-q_0|^2-Mt.
\end{align*}
If $t<\sigma$, then the equality 
$
S^t(q_0,q_1)=A^t(q_0,q_1;n)
$
can be proved by recurrence for each $n$ using Proposition \ref{triangle}.
For general $t$, let us prove that $A^t(n)$ is independent of $n$. We take two integers $n$ and $m$ such that 
$t/n<\sigma, t/m<\sigma$ and want to prove that 
$A^t(n)=A^t(m)$. We will prove that 
$A^t(n)=A^t({nm})=A^t(m)$.
Since $t/m<\sigma$, we have
$$
A_{\tau}^{\tau+t/m}(q_0,q_1;n)=
S_{\tau}^{\tau+t/m}(q_0,q_1)
$$
for each $\tau$ and $n$, hence
\begin{align*}
A^t(q_0,q_1;nm)=\quad&\min _{\theta_1,\theta_2,\ldots \theta_{nm-1}}
\Big[S_0^{t/nm}(q_0,\theta_1)+S_{t/nm}^{2t/mn}(\theta_1,\theta_2)+\cdots+
S_{(n-1)t/nm}^{t/m}(\theta_{n-1},\theta_n)\\
&
+S_{t/m}^{(n+1)t/nm}(\theta_n,\theta_{n+1})+\cdots+
S_{(2n-1)t/nm}^{2t/m}(\theta_{2n-1},\theta_{2n})\\
& +\cdots\\
&
+S_{(m-1)t/m}^{(m-1)t/m+t/nm}(\theta_{(m-1)n},\theta_{(m-1)n+1})+\cdots+
S_{(1-1/nm)t}^{t}(\theta_{mn-1},q_1)
\Big]\\
=\quad &\min_{\theta_{2n},\theta_{3n},\ldots ,\theta_{(m-1)n}}
\Big[S_0^{t/m}(q_0,\theta_n)+S_{t/m}^{2t/m}(\theta_n,\theta_{2n})+\cdots+
S_{(m-1)t/m}^{t}(\theta_{(m-1)n},q_1)\Big]\\
=\quad & A^t(q_0,q_1;m).
\end{align*}
We have proved that $A^t({nm})=A^t(m)$, by symmetry we also have
$A^t({nm})=A^t(n)$ hence $A^t(n)=A^t(m)$.
Finally, we have 
$$
\mathfrak{S}(\theta_1,\ldots,\theta_{n-1})
\leq \frac{n}{2mt}
\big(|\theta_1-q_0|^2+|\theta_2-\theta_1|^2
+|q_1-\theta_{n-1}|^2\big)+Mt
$$
hence 
\begin{align*}
A^t(q_0,q_1)&\leq \min _{\theta_1,\theta_2,\ldots \theta_{n-1}}\frac{n}{2mt}
\big(|\theta_1-q_0|^2+|\theta_2-\theta_1|^2
+|q_1-\theta_{n-1}|^2\big)+Mt\\
&=
\frac{1}{2mt}|q_1-q_0|^2+Mt.
\end{align*}
\qed
The following property concerning $A$ follows easily from
the definition :
\begin{equation}\label{triangleA}\tag{T}
A_{t_0}^{t_2}(q_0,q_2)=\min_{q_1} \big( 
A_{t_0}^{t_1}(q_0,q_1)+A_{t_1}^{t_2}(q_1,q_2)
\big),
\end{equation}
when $0\leq t_0\leq t_1\leq t_2$.
The following consequence of Hypothesis \ref{h3}
 will also be useful:
\begin{lem}\label{estimeeaction}
$$p\cdot \partial_pH(t,q,p)-H(t,q,p)\geq \frac{m}{M}H(t,q,p)
-(m+M).$$
\end{lem}
\proof
We deduce from Hypothesis \ref{h2} that 
$$
H(t,q,0)\geq H(t,q,p)-p\cdot \partial_pH(t,q,p)+\frac{m}{2}|p|^2.
$$
We deduce that 
$$
p\cdot \partial_pH(t,q,p)-H(t,q,p)\geq \frac{m}{2}|p|^2
-H(t,q,0)\geq \frac{m}{M}(H(t,q,p)-M)
-M
$$
\qed

The minimal action $A^t(q_0,q_1)$ is not necessarily
$C^1$, we need some definitions before we can study
its regularity.
The linear form $l$ is called a $K$-super-differential 
of the function $u$ at point $q$ if 
the inequality
$$
u(\theta)\leq u(q)+l(\theta-q)+K|\theta-q|^2
$$ 
holds in a neighborhood of $q$.
The linear form $l$ is a proximal super-differential of $u$ at point
$q$
if it is a $K$-super-differential for some $K$.
The form $l$ is a proximal super-differential of $u$ at $q$
if and only if there exists a $C^2$ function $v$
such that $dv(q)=l$ and such that the difference $v-u$ 
has a minimum at $q$.
More generally, we will say that $l$ is 
a super-differential of $u$ at $q$ if 
there exists a
$C^1$ function $v$
such that $dv(q)=l$ and such that the difference $v-u$ 
has a minimum at $q$.
A super-differential is not necessarily a proximal super-differential.

A function $u:\Rm^d\lto \Rm$  is called $K$-semi-concave if
it admits a $K$-super-differential at each point.
It is equivalent to require that the function 
$\theta\lmto u(\theta) -K|\theta|^2$ is concave.
A function is called semi-concave if it is $K$-semi-concave for some $K$.
If $u$ is a $K$-semi-concave function,
and if $l$ is a super-differential at $u$,
then the inequality
$$
u(\theta)\leq u(q)+l(\theta-q)+K|\theta-q|^2
$$ 
holds for each $\theta$. In particular,
$l$ is a $K$-super-differential.

\begin{lem}
The function 
$A^t$
is $C(1+1/t)$-semi-concave, with some constant $C$
which depends only on $m$ and $M$.
\end{lem}

\proof
Let us first assume that $t\in ]0,\sigma[$.
In this case, $A_0^t=S_0^t$, this function is $C^2$ and its second derivative was estimated in Lemma \ref{estimeesS}.
Let us now assume that $t\geq \sigma$.
Then, there exists $n\in \Nm$ such that $t/n\in[\sigma/3,\sigma/2[$.
We  have 
$$
A_0^t(q,q')=
\min _{\theta,\theta'} \big(S_0^{t/n}(q,\theta)+
A_{t/n}^{t-t/n}(\theta,\theta')
+S_{t-t/n}^{t}(\theta',q')
\big).
$$
Considering a minimizing pair $(\theta_0, \theta_1)$
in the expression above at $(q_0,q_1)$, we see that 
the $C^2$ function 
$$
(q,q')\lmto S_0^{t/n}(q,\theta_0)+A_{t/n}^{t-t/n}(\theta_0,\theta_1)
+S_{t-t/n}^{t}(\theta_1,q')
$$
is touching from above the function $A_0^t$ at point $(q_0,q_1)$.
In view of Lemma \ref{estimeesS},
this provides a uniform (for $t\geq \sigma$) semi-concavity constant 
for $A_0^t$.
\qed

\section{The Lax-Oleinik operators.}
Given $t_0<t_1$, we define the Lax-Oleinik
operators ${\mathbf{T}}_{t_0}^{t_1}$ and
$\check {\mathbf{T}}_{t_1}^{t_0}$
which, to each function $u:\Rm^d\lto \Rm$
associate the functions 
$$
{\mathbf{T}}_{t_0}^{t_1}u(q):= \inf_{\theta \in \Rm^d} \big( u(\theta)+
A_{t_0}^{t_1}(\theta,q)\big)
\quad ,\quad
\check {\mathbf{T}}_{t_1}^{t_0}u(q):= \sup_{\theta \in \Rm^d} \big( u(\theta)-
A_{t_0}^{t_1}(q,\theta)\big).
$$
We have the Markov (or semi-group) property:
$$
{\mathbf{T}}^{t_2}_{t_1}\circ {\mathbf{T}}^{t_1}_{t_0}={\mathbf{T}}_{t_0}^{t_2}
\quad,\quad
\check {\mathbf{T}}^{t_0}_{t_1}\circ\check {\mathbf{T}}^{t_1}_{t_2}=
\check {\mathbf{T}}_{t_2}^{t_0}
$$
for $t_0<t_1<t_2$.
Note however that 
${\mathbf{T}}^{t_1}_{t_0}\circ \check {\mathbf{T}}^{t_0}_{t_1}$ 
and 
$
\check 
{\mathbf{T}}^{t_0}_{t_1}\circ
{\mathbf{T}}^{t_1}_{t_0}
$
are not the identity.
Concerning these operators, we only have the inequalities
$$
\check 
{\mathbf{T}}^{t_0}_{t_1}\circ
{\mathbf{T}}^{t_1}_{t_0}(u)
\leq u
\quad,\quad
{\mathbf{T}}^{t_1}_{t_0}\circ \check {\mathbf{T}}^{t_0}_{t_1} (u)\geq u,
$$
the easy proof of which is left to the reader.
Each property concerning the Lax-Oleinik operator
${\mathbf{T}}$  has a 
counterpart for the dual operator
$\check {\mathbf{T}}$, that we will not always  bother to state,
but never  hesitate to use.
The family of operators ${\mathbf{T}}_{t_0}^{t_1}$ is characterized by
the fact that
 ${\mathbf{T}}_{t_0}^{t_1}u(q)=\inf_{\theta} \big(u(\theta)+
S_{t_0}^{t_1}(q_0,q_1)\big)
$ 
when $t_0\leq t_1\leq t_0+\sigma$
and by the Markov property.
The Lax-Oleinik operators solve (\ref{HJ})
in various important ways, that will be detailed in the present section. 
It is useful first to settle some regularity issues.

\begin{lem}\label{semiconst}
There exists a constant $C$, depending only on $m$
and $M$ such that, for each 
$t\in ]0,\sigma]$, the function  ${\mathbf{T}}^tu$ is
$(C/t)$-semi-concave provided it is finite at each point.
\end{lem}

\proof
The function ${\mathbf{T}}^tu$ is the infimum of the functions
$f=u(\theta)+S^t(\theta,.)$, which are $C^2$ with
the uniform bound $\|d^2f\|\leq C/t$.
It is then an easy exercise to conclude that the function 
${\mathbf{T}}^tu$ is $C/t$-semi-concave, see Lemma \ref{infimum}.
\qed

Given an arbitrary function $u_0$, the infimum in 
the definition of ${\mathbf{T}}_0^t u_0$ is not necessarily
finite, and, even if it is finite, it is not necessarily
a minimum.
It is clear from Proposition \ref{estimeA}
that the infimum is a finite minimum
under the assumption that $u_0$ is continuous and
\textbf{Lipschitz in the Large}, which means that
there exists a constant $k$ such that 
$$
u_0(q')-u_0(q)\leq k(1+|q'-q|)
$$
for each $q$ and $q'$.

\begin{lem}\label{lil}
If $u_0$ is Lipschitz in the large, then 
so are the functions ${\mathbf{T}}^t_0u_0$ for all
 $t\geq 0$.
The function $(t,q)\lmto u(t,q)= {\mathbf{T}}_0^tu_0(q)$
is locally semi-concave, 
hence locally Lipschitz on $]0,\infty)\times \Rm^d$.
The function $u$ solves
(\ref{HJ}) at all its points of differentiability (hence almost everywhere).
\end{lem}

\proof
Since $u_0$ is Lipschitz in the large, 
the function ${\mathbf{T}}^t_0u_0-u_0$ is bounded
for each $t>0$,
as follows from the inequalities
$$
\inf_{\theta}
\big(
u_0(q)-k-k|\theta-q|+S^t(\theta,q)
\big)
\leq {\mathbf{T}}^t_0u_0 \leq u_0(q)+S^t(q,q)
$$
which imply (setting $\Delta=\theta-q$) that 
$$
\inf_{\Delta\in \Rm^d}
\big(
-k-k|\Delta|+\frac{1}{2tM}|\Delta|^2-tM
\big)
\leq
{\mathbf{T}}^t_0u_0(q)-u_0(q)
\leq Mt.
$$
We conclude that the function ${\mathbf{T}}^t_0u_0=({\mathbf{T}}_0^tu_0-u_0)+u_0$
is
Lipschitz in the large.
In the computations above, we also see that 
the infimum can be taken on $|\Delta|\leq K$,
where $K$ is a constant independent from $q$.

Let us now prove that the function $u(t,q):={\mathbf{T}}_0^tu_0(q)$
is locally Lipschitz on $t>0$.
 In view of the Markov property,
it is enough to prove that the function $u$ is 
Lipschitz on $]\tau,\sigma/2[\times B$
for each closed ball $B\subset \Rm^d$ and each  time
$\tau\in ]0,\sigma/2[$.
Since $u(q)$ is Lipschitz in the large,
there exists a radius  $R>0$ such that 
$$
u(t,q)=\inf_{|\theta|\leq R} u(\theta)+S^t(\theta,q)
$$
for $(t,q)\in ]\tau,\sigma/2[\times B$.
Since $S$ is $C^2$, the functions 
$(t,q) \lmto u(\theta) + S^t(\theta,q),|\theta|\leq R$
have uniform $C^2$ bounds  on $]\tau,\sigma/2[\times B$.
Their infimum $u(t,q)$ is then semi-concave,
hence Lipschitz  on that set, see Lemma \ref{infimum}.

Finally, let $(t,q)$ be a point of differentiability of $u$,
and let $\tau \in ]\max(0,t-\sigma),t[$ be given.
Since $u_{\tau}$ is Lipschitz in the large and locally Lipschitz,
there exists $\theta$ such that 
$\mathbf{T}_{\tau}^t u_{\tau}(q)=u_{\tau}(\theta)+S_{\tau}^t(\theta,q)$.
For a different point $(s, y)$, we have 
$\mathbf{T}_{\tau}^s u_{\tau}(y)\leq u_{\tau}(\theta)+S_{\tau}^t(\theta,y)$,
hence the function $(s,y)\lmto u(s,y)-S_{\tau}^s(\theta,y)$
has a maximum at $(t,q)$, which implies that the functions $u(s,y)$ and 
$S^s_{\tau}(\theta,y)$, each of which is differentiable at $(t,q)$, 
 have the same differential at $(t,q)$.
Since the functions $(s,y)\lmto S_{\tau}^s(\theta,y)$ 
solves (\ref{HJ}), the function $u$ also solves (\ref{HJ}) at $(t,q)$.
\qed

Let us  now establish the relation of our operators with
regular solutions.

\begin{prop}\label{formel}
Let  $u(t,q):]t_0,t_1[\times \Rm^d \lto \Rm$ be a  $C^1$ solution
of \ref{HJ},
 then 
${\mathbf{T}}_{\tau}^tu_{\tau}=u_t$ 
and $\check {\mathbf{T}}^{\tau}_tu_{t}=u_{\tau}$ 
for each $\tau \leq t$ in $]t_0,t_1[$.
The function $u$ is locally $C^{1,1}$.
 \end{prop}
This property is one of the main motivations to introduce the Lax-Oleinik operators.
The observation that $C^1$ solutions
are actually locally $C^{1,1}$  comes Fathi's paper 
\cite{Fa:03}, itself inspired by anterior works of  Herman.
Another consequence of this Theorem 
is that uniqueness extends to $C^1$ solutions
under the convexity assumption.

\proof
In view of the Markov property, 
it is enough to prove the result for $0<t-\tau<\sigma$.
Given $q$ and $\theta$ in $\Tm^d$, we consider the unique
orbit $(q(s),p(s))$ such that $q(\tau)=\theta$
and $q(t)=q$.
By the convexity of $H$, we have
$$
H(q(s),\partial_q u(s,q(s)))\geq H(q(s),p(s))+
(\partial_qu(s,q(s))-p(s))\cdot \partial_pH(s,q(s),p(s)).
$$
Noticing that $\dot q(s)=\partial_pH(s,q(s),p(s))$
and integrating gives:
\begin{align*}
S_{\tau}^t(\theta,q)&=\int_{\tau}^t p(s)\cdot \dot q(s) -H(s,q(s),p(s)) ds\\
&\geq 
\int_{\tau}^t \partial_q u(s,q(s))\cdot \dot q(s) -H(s,q(s), \partial_q u(s,q(s))) ds\\
&=u(t,q)-u(\tau,\theta),
\end{align*}
with equality if $p(s)= \partial_q u(s,q(s))$ for each $s$.
We conclude that
$$
{\mathbf{T}}_{\tau}^t u_{\tau}(q)\geq u_t(q),
$$
with equality if there exists an orbit 
$(q(s),p(s)): [\tau,t]\lto \Rm^d\times \Rm^{d*}$
such that $p(s)= \partial_q u(s,q(s))$ and $q(t)=q$.
By Corollary \ref{geomc1}, the orbit of the point
$(q,\partial_q u(t,q))$ satisfies this property,
hence the equality holds.

To prove the regularity of $u$ we consider
a subinterval $[\tilde t_0,\tilde t_1]\subset ]t_0,t_1[$,
and prove that $u$ is locally $C^{1,1}$ on $]\tilde t_0,\tilde t_1[$.
We have 
$$
u(t,q)={\mathbf{T}}_{\tilde t_0}^t u_{\tilde t_0}(q)=\
\check {\mathbf{T}}^{\tilde t_1}_t u_{\tilde t_1}(q)
$$
for each $t\in ]\tilde t_0,\tilde t_1[$.
If the functions $u_t$ were Lipschitz in the Large, we could
apply Lemma \ref{lil} and deduce that $u$ is both 
locally semi-concave and locally semi-convex, hence locally $C^{1,1}$,
 on  $]\tilde t_0,\tilde t_1[\times \Rm^d$.
Here we do not make any growth assumption, so we need a slightly different argument
to prove the semi-concavity of $u$ (and, similarly, its semi-convexity).
We have seen that 
the infimum in the definition ${\mathbf{T}}_{\tilde t_0}^t u_{\tilde t_0}(q)$
is a minimum, which is attained at the point 
$\theta = Q^{\tilde t_0}_t (q, \partial_q u(t,q)).
$
This gives us an \textit{a priori} bound on $\theta$, and we can continue the proof
as in Lemma \ref{lil}.
\qed

Let us sum up some  properties of the 
Lax-Oleinik operators ${\mathbf{T}}_{\tau}^t$
associated to a Hamiltonian
satisfying hypotheses 
\ref{h1},\ref{h2},\ref{h3}: 

\begin{property}\label{axiom} $ $
\begin{enumerate}
\item Markov property:
 ${\mathbf{T}}_s^t\circ {\mathbf{T}}_{\tau}^s={\mathbf{T}}_{\tau}^t$
 when $\tau\leq s\leq t$.
\item Monotony: $u\geq v  \Rightarrow
 {\mathbf{T}}_{\tau}^tu\geq {\mathbf{T}}_{\tau}^tv$ for each $t\geq \tau$.
\item Compatibility with (\ref{HJ}): 
If $u(t,q):]t_0,t_1[\times \Rm^d\lto \Rm$ is a $C^2$
solution of (\ref{HJ}), then ${\mathbf{T}}_{\tau}^tu_{\tau}=u_t$ when
$t_0<\tau<t<t_1$.
\item Boundedness: If $u_{\tau}$ is Lipschitz in the large,
then the functions ${\mathbf{T}}_{\tau}^tu_{\tau}$,
$t\in [\tau,T]$ are uniformly Lipschitz in the large
for each $T\geq \tau$.
\item Regularity:  If $u_{\tau}$ is Lipschitz in the large, the function 
$(t,q)\lmto {\mathbf{T}}_{\tau}^tu_{\tau}(q)$ is locally Lipschitz
on $]\tau,\infty)\times \Rm^d$.
\item Translation invariance: 
${\mathbf{T}}_{\tau}^t(c+u)=c+{\mathbf{T}}_{\tau}^tu$ for each constant $c\in \Rm$.
 \end{enumerate}
\end{property}

The Lax-Oleinik operators solve the Cauchy problem for 
(\ref{HJ}) in the viscosity sense.
Actually, this  follows from Property \ref{axiom}:

\begin{prop}\label{viscous}
Let $H$ be a Hamiltonian satisfying Hypothesis \ref{h1}.
Assume that there exists a family ${\mathbf{T}}_{\tau}^t$, 
$0\leq \tau\leq t$ 
of operators satisfying the Markov property, the monotony, the compatibility with (\ref{HJ}), and the boundedness as expressed in Property \ref{axiom}.
Then if  $u_{0}$ is an initial condition which is Lipschitz in
the large, 
the function
$$
(t,q)\lmto u(t,q)={\mathbf{T}}_{0}^tu_{0}(q)
$$
is a viscosity solution of (\ref{HJ}) on 
$]0,\infty)\times \Rm^d$.
\end{prop}

Notice that we did not make any convexity assumption.
This kind of axiomatic characterization of viscosity
solutions is reminiscent from \cite{AGLM}, see also
\cite{Bit}. It may also help to understand the links between viscosity solutions and variational
solutions in the non-convex setting. Such links were suggested 
by Claude Viterbo and established in her thesis
by Qiaolin Wei, \cite{QW}.

\textsc{Proof of Proposition \ref{viscous}: }
Let us prove that $u$ is a viscosity sub-solution,
a similar proof yields that it is also a super-solution.
We consider a point  $(T, Q)\in ]0,\infty)\times \Rm^d$
and  a super-differential $(h,p)$   
of the function $u$ at $(T,Q)$.
To prove that $h+H(T,Q,p)\leq 0$,
we assume, by contradiction,  that 
$$h+H(T,Q,p)>0.
$$
As is usual for viscosity solutions
we will use a test function $\phi$.
We will assume that $\phi:\Rm\times \Rm^d\lto \Rm$ is smooth and satisfies  the following properties:
\begin{itemize}
\item
$
\phi(T,Q)=u(T,Q),\quad 
 \partial_t\phi(T,Q)=h,\quad
\partial_q\phi(T,Q)=p,
$
\item $\phi\geq u$ on $[-T/2,2T] \times \Rm^d$,
\item There exists a constant $C>0$ such that 
$\phi(t,q)=C\sqrt{1+|q|^2}$ when $|q|+|t|\geq C$.
\end{itemize}
Note that $d^2\phi$ is bounded.  Such a test function exists because the functions 
$u_{t}$, $t\in [T/2,2T]$, are uniformly Lipschitz in the 
large, as follows from the boundedness property assumed
on the operators.

\noindent
\textbf{Claim:} There exists  $S>0$ and a $C^2$ function 
$w(\tau,t,q)$ defined on the open set
$$
\big\{
(\tau,t,q)\in\Rm\times \Rm \times \Rm^d: \tau-S<t<\tau+S
\big\}
\subset\Rm\times \Rm \times \Rm^d
$$
such that, for each fixed $\tau$, the function 
$w_{\tau}:(t,q)\lmto w(\tau,t,q)$ is the solution
of the Cauchy problem
$$
\left\{
\begin{matrix}
\partial_t w_{\tau} 
+H(t,q,\partial_q w_{\tau})=0
\\
w_{\tau}(\tau,q)=\phi(\tau,q).
\end{matrix}
\right.
$$
The existence of a solution $w_{\tau}$
to this problem  follows from Theorem \ref{existence}.
However, to see that $w$ is $C^2$ in all its variables,
we find it  more convenient
to consider  the Cauchy problem
$$
\left\{
\begin{matrix}
\partial_s u +
\Big(\partial_ z u
+H(z,q,\partial_q u(s,z,q))
\Big)=0
\\
u(0,z,q)=\phi(z,q).
\end{matrix}
\right.
$$
By Theorem  \ref{existence},  
applied to the Hamiltonian 
\begin{align*}
\hat H(s,z,q,\xi,p):\Rm\times 
(\Rm \times  \Rm^d)\times (\Rm\times \Rm^{d})^*
&\lto \Rm\\
(s,z,q,\xi,p)&\lmto \xi+H(z,q,p)
\end{align*}
there exists $S>0$
and a $C^2$ solution 
$u(s,z,q):]-S,S[\times \Rm\times \Rm^d\lto \Rm$
of this Cauchy problem.
Setting 
$$
w(\tau,t,q):= u(t-\tau,t,q),
$$
we verify that 
$$
\partial_tw(t,q)+H(t,q,\partial_qw(t,q))=
\partial_su(t-\tau,t,q)+\partial_z u(t-\tau,t,q)
+H(t,q,\partial_qu(t-\tau,t,q))=0
$$
and that 
$
w(\tau,\tau,q)=u(0,\tau,q)=\phi(\tau,q).
$

\noindent
\textbf{Claim :} There exists $\tau\in ]T-S,T[$
such that
$w(\tau,T,Q)<\phi(T,Q)$. 

Since $w(T,T,q)=\phi(T,q)$, we have 
$$
\partial_t w(T,T,Q)=-H(T,Q,\partial_qw(T,Q))
=-H(T,Q,\partial_q\phi(T,Q))<\partial_t\phi(T,Q).
$$
As a consequence, there exists $\delta>0$ such that
$$
\partial_tw(\tau,t,Q)-\partial_t\phi(t,Q)<0
$$
for $\tau,t\in ]T-\delta,T[$.
Since $w(\tau,\tau,Q)=\phi(\tau,Q)$, we deduce
by integration  that 
$$
w(\tau,T,Q)-\phi(T,Q)=
\int_{\tau}^T \partial_tw(\tau,t,Q)-
\partial_t\phi(t,Q)dt
<0
$$
provided $\tau\in ]T-\delta,T[$, which proves our claim.

\noindent
\textbf{Conclusion :}
Since we are considering monotone  operators compatible
with (\ref{HJ}) we have 
$$
w(\tau,T,Q)={\mathbf{T}}_{\tau}^Tw_{\tau}
(Q)={\mathbf{T}}_{\tau}^T\phi_{\tau}(Q)\geq {\mathbf{T}}_{\tau}^Tu_{\tau}(Q)=u(T,Q) 
$$
hence $\phi(T,Q)>u(T,Q)$, which is a contradiction.
\qed

This parenthesis through viscosity solutions
being closed,
let us  turn our attention to more geometric 
aspects of the Lax-Oleinik operators.
 We denote 
by $\Gamma_u$ the graph of the differential of $u$ on its
domain of definition,
$$
\Gamma_u:= \{(q,du(q)): \quad q\in \Rm^d, \quad  du(q) \,\text{exists}\}.
$$

\begin{prop}\label{geom}
Let $u$ be a semi-concave and Lipschitz function.
The set 
$$
\varphi_t^0\left(\bar\Gamma_{{\mathbf{T}}_0^tu}
\right)
$$
is contained in $\Gamma_u$ for each $t>0$, and it is a Lipschitz graph.
\end{prop}

\proof
In view of the Markov property, it is enough to prove the result for 
$t\in ]0,\sigma]$.
Let $(q,p)$ be a point of $\Gamma_{{\mathbf{T}}^t_0 u}$, which means
that the function ${\mathbf{T}}_0^t u$ is differentiable at $q$
and that $d({\mathbf{T}}_0^t u)(q)=p$.
Let $\Theta$ be a minimizing point in the expression
${\mathbf{T}}_0^tu(q)=\min_{\theta} u(\theta)+S_0^t(\theta,q)$.
Since each of the functions $u$ and $S_0^t(.,q)$
are semi-concave, this implies that they are both differentiable at $\Theta$, and that
 $du(\Theta)+\partial_0 S_0^t(\Theta,q)=0$.
Moreover, this implies that the function
$ u(\Theta)+S_0^t(\Theta,.)
$
touches the function ${\mathbf{T}}_0^tu$ from above at point $q$,
hence that $S_0^t(\Theta,.)$
is differentiable at $q$, with a differential equal to $p$.
We then have
$$
\varphi_t^0(q,p)=\varphi^0_t(q,\partial_1S_0^t(\Theta,q))
=(\Theta,- \partial_0S_0^t(\Theta,q)=(\Theta,du(\Theta))
\subset \Gamma_u.
$$
We have proved that $\varphi_t^0(\Gamma_{{\mathbf{T}}_0^tu})\subset \Gamma_u$.
Moreover, we have 
$Q_t^0(\Gamma_{{\mathbf{T}}_0^tu})\subset \mI$, where $\mI\subset \Rm^d$
is the set of points $\theta$ which are minimizing in the definition 
of ${\mathbf{T}}_0^tu(q)$ for some point $q$.

\noindent
\textbf{Claim: } The function  $u$ is $C^{1,1}$ on $\mI$. This  means that
$u$ is differentiable at each point of $\mI$, and that the 
map $\theta\lmto du(\theta)$ is Lipschitz on $\mI$.
In other words, the projection of $\Gamma_u$ to $\Rm^d$
contains $\mI$, and the set 
$$
\Gamma_{u|\mI}:=\{(\theta, du(\theta), \theta\in \mI\}
$$ 
is a Lipschitz graph.

To prove the claim, we first prove that $u$ has $C$-super-differentials
and $C$-sub-differentials at each point of $\mI$, where $C$ is a 
common semi-concavity constant of all the functions 
$-S_0^t(.,q)$ and of the function $u$.
The existence of a $C$-super-differential follows from the $C$-semi-concavity of $u$.
To prove the existence of a $C$-sub-differential at a point $\Theta\in \mI$,
we consider a point $q$ such that $u(\Theta)+S_0^t(\Theta,q)=T_0^tu(q)$.
Such a point exists by definition of $\mI$.
This implies that the function $\theta\lmto u(\theta)+S_0^t(\theta,q)$ has a minimum 
at $\theta=\Theta$, hence  each $C$-sub-differential of $-S_0^t(.,q)$ is
 a $C$-sub-differential 
of $u$.
The claim then follows from a result of Fathi, see Proposition
 \ref{fathi}
in the Appendix.

Let now  $(q,p)$ be a point in the closure 
$\bar\Gamma_{{\mathbf{T}}_0^tu}$ of $\Gamma_{{\mathbf{T}}_0^tu}$.
There exists a sequence $(q_n,p_n)$ of points of 
$\Gamma_{{\mathbf{T}}_0^tu}$ which converges to $(q,p)$.
By definition, the function ${\mathbf{T}}_0^tu$ is differentiable at $q_n$,
and $p_n=d({\mathbf{T}}_0^tu)(q_n)$.
Let $\Theta_n=Q_t^0(q_n,p_n)$ be the sequence of points such that 
$$
{\mathbf{T}}_0^tu(q_n)=u(\Theta_n)+S_0^t(\Theta_n,q_n).
$$
The sequence $\Theta_n$ is converging to $\Theta=Q_t^0(q,p)$,
and, at the limit, we see that
$$
{\mathbf{T}}_0^tu(q)=u(\Theta)+S_0^t(\Theta,q).
$$
We conclude that $\Theta\in \mI$.
Since we have already proved the Lipschitz regularity
of $du$ on $\mI$, we deduce that 
$
\varphi_t^0(q,p)=\lim (\varphi_t^0(q_n,p_n))=\lim (\Theta_n,du(\Theta_n))=
(\Theta,du(\Theta))\in
\Gamma_{u|\mI}\subset \Gamma_u.
$
\qed

The action of the Lax-Oleinik operators on semi-convex functions
also has a remarkable property, see \cite{ENS}.
It is useful to denote by $L_u$ the set 
of point 
$(Q,P)$ such that 
$P$ is a  sub-differential of $u$ at $Q$.
Note that $\Gamma_u\subset L_u$.
\begin{prop}\label{semiconvex}
If $u$ is $K$-semi-convex, then for
each $\delta\in]0,1[$ there exists $T>0$ such that ${\mathbf{T}}^t_0u$
is $(K+\delta)$-semi-convex, hence $C^{1,1}$, for each $t\in ]0,T]$.
One can take 
$$
T=\frac{\delta }{2M(3+2K)^2}.
$$
\end{prop}
\proof
Since $u$ is $K$-semi-convex, for each $(Q,P)\in L_u$, we have 
$$
u(q)\geq u(Q)+P(q-Q)-K|q-Q|^2.
$$
We denote by $l_{Q,P}(q)$ the function on the right in this inequality, so that 
$$
u=\max_{(Q,P)\in L_u} l_{Q,P}.
$$
Taking $T$ as in the statement, it follows from Theorem \ref{existence} that
the functions ${\mathbf{T}}^t_0(l_{Q,P}), t\in [-T,T]$
 are $C^2$ with a second derivative bounded by $2K+4tM(1+2K)^2\leq 2K+2\delta$.
We claim that 
$$
{\mathbf{T}}^t_0u=\max_{(Q,P)\in L}{\mathbf{T}}^t_0(l_{Q,P}),
$$
for $t\in[0,T]\cap[0, \sigma]$,
which implies that ${\mathbf{T}}^t_0u$ is $(K+4tM(1+2K)^2)$-semi-convex.
We prove the claim in two steps. First, the inequality
$$
{\mathbf{T}}^t_0u\geq\max_{(Q,P)\in L}{\mathbf{T}}^t_0(l_{Q,P})
$$
follows immediately from the fact that $u\geq l_{Q,P}$ for each $(Q,P)\in L$
in view of the monotony of ${\mathbf{T}}^t_0$, see Property \ref{axiom}.
Let us fix a point $(t,q)$ and prove the converse inequality at this point.
Since 
$$
u(\theta)+S^t_0(\theta,q)\geq u(q)+P(\theta-q)-K(\theta-q)^2+
\frac{1}{2tM}|\theta-q|^2
-tM
$$
and since $K\leq 1/2tM$, there exists a point $\theta$ such that
${\mathbf{T}}^t_0u(q)=u(\theta)+S^t_0(\theta,q)$.
Assuming that $t\leq \sigma$, 
this implies that the point 
$(\theta,\zeta)=(\theta, -\partial_0S^t_0(\theta,q))$ belongs to $L_u$,
and that $q=Q^t_0(\theta,\zeta)$.
Then, we have 
$${\mathbf{T}}^t_0(l_{\theta,\zeta})(q)=l_{\theta,\zeta}(\theta)+S^t_0(\theta,q)
=u(\theta)+S^t_0(\theta,q)={\mathbf{T}}^t_0u(q)
$$
hence 
$${\mathbf{T}}^t_0u(q)\leq\max_{(Q,P)\in L} {\mathbf{T}}^t_0(l_{Q,P})(q)
$$
provided $t\leq \sigma$.
We conclude that ${\mathbf{T}}^t_0u$ is semi-concave with constant
$K+2tM(1+2K)^2$ for $t\in[0,\sigma]\cap [0,T]$.
We can then apply this result to ${\mathbf{T}}_0^{\sigma}u$, and,
since $K+tM(1+2K)^2\leq K+1$, we conclude that
the function
${\mathbf{T}}_{\sigma}^t{\mathbf{T}}_0^{\sigma}u$ is semi-concave with  constant
$$
K+2\sigma M(1+2K)^2+2tM(3+2K)^2\leq K+2(\sigma+t)M(3+2K)^2\leq K+1
$$
 for 
$t\in [0, \sigma]\cap [0, T-\sigma]$.
In other words, the functions 
${\mathbf{T}}_0^tu$ are semi-concave with constant 
$K+2tM(3+2K)^2$ for $t\in [0,2\sigma]\cap [0,T]$.
We can apply this argument as many times as necessary
and obtain that, 
the functions ${\mathbf{T}}_0^{t}u$ are semi-concave with constant
$K+2tM(3+2K)^2$ for each $t\in [0,T]$.
\qed
The following was first stated explicitly by Marie-Claude Arnaud in \cite{Ar:11}.

\begin{ad}
Under the hypotheses of Proposition \ref{semiconvex},
we have $L_u=\varphi^{0}_t(\Gamma_{{\mathbf{T}}^t_0u})$ for each $t\in]0,T[$.
Moreover, for each $q$, we have ${\mathbf{T}}_0^tu(q)=u(\theta)+
S^t_0(\theta,q)$, with $\theta=Q_t^{0}(q,d({\mathbf{T}}^t_0u)(q))$.
\end{ad}
\proof
For each $q\in \Rm^d$, we have seen that there exists 
$(\theta,\zeta)\in L_u$ such that 
${\mathbf{T}}^t_0u(q)=u(\theta)+S^t_0(\theta, q)$
and $\zeta=-\partial_0S^t_0(\theta,q)$.
Since we know that ${\mathbf{T}}_0^tu$ is $C^1$, the first of these equalities
implies that $d({\mathbf{T}}^t_0u)(q)=\partial_1S^t_0(\theta,q)$, while the second implies
that 
$\varphi^t_0(\theta,\zeta)=(q,\partial_1S^t_0(\theta,q))$.
We conclude that $\varphi^{0}_t(\Gamma_{{\mathbf{T}}^t_0u})\subset L_u$.
Moreover,  $\theta=Q^{0}_t(q,d({\mathbf{T}}^t_0u)(q))$.

Conversely, let us consider a point $(\theta,\zeta)\in L$, and 
denote by $l$ the associated
function $l_{\theta,\zeta}$.
By Proposition \ref{formel},
 the function $(t,q)\lmto {\mathbf{T}}^t_0l(q)$ is the restriction to $]0,T[\times \Rm^d$
of the  $C^2$
solution of (\ref{HJ}) emanating from $l$.
As a consequence, we have
$$
{\mathbf{T}}^t_0l(Q^t_0(\theta,\zeta))=l(\theta)+S^t_0(\theta,Q^t_0(\theta,\zeta))
=u(\theta)+S^t_0(\theta,Q^t_0(\theta,\zeta))\geq {\mathbf{T}}^t_0 u(Q^t_0(\theta,\zeta)).
$$
Since we know from the monotony property that ${\mathbf{T}}^t_0l\leq {\mathbf{T}}^t_0u$, we conclude that this last inequality is actually an equality.
Setting $q_1=Q^t_0(\theta,\zeta)$, this  implies that 
$$
(\theta, \zeta)=(\theta,-\partial_0S_0^t(\theta, q_1))=
\varphi^{0}_t(q_1,\partial_1S^t_0(\theta,q_1))=
\varphi^{0}_t(q_1,d{\mathbf{T}}^t_0u(q_1))
\subset \varphi^{0}_t(\Gamma_{{\mathbf{T}}^t_0u}).
$$
We conclude that $L_u\subset  \varphi^{0}_t(\Gamma_{{\mathbf{T}}^t_0u})
.$
\qed

\begin{ad}\label{ad}
Under the hypotheses of Proposition \ref{semiconvex},
we have $\check {\mathbf{T}}_t^0 \circ {\mathbf{T}}_0^tu=u$ for each $t\in ]0,T[$.
\end{ad}

\proof
Let us define the map $F:q\lmto Q_t^0(q,d({\mathbf{T}}_0^tu(q))$.
By the first addendum, the image of $F$ is equal to the projection
of $L_u$ on $\Rm^d$, hence the map $F$ is onto.
Given a point $\theta\in \Rm^d$,
we consider a preimage $q$ of $\theta$ by $F$,
and write
$$
\check {\mathbf{T}}_t^0 \circ {\mathbf{T}}_0^t u (\theta)\geq 
{\mathbf{T}}_0^tu(q)-S_0^t(\theta, q)
=
u(\theta)
$$ 
where the last equality comes from the first addendum.
We conclude that
$\check {\mathbf{T}}_t^0 \circ {\mathbf{T}}_0^t u\geq u$, hence 
$\check {\mathbf{T}}_t^0 \circ {\mathbf{T}}_0^t u= u$.
\qed

The following extrapolates on \cite{Padova}.
For $t_0\in \Rm$ and $\delta, t>0$, let us define the operators
$$
{\mathbf{R}}^t:=\check {\mathbf{T}}_{t_0+\delta t}^{t_0}
\circ {\mathbf{T}}_{t_0-t}^{t_0+\delta t}
\circ \check {\mathbf{T}}_{t_0}^{t_0-t}
\quad , \quad
\check {\mathbf{R}}^t :={\mathbf{T}}^{t_0}_{t_0-\delta t}
\circ \check {\mathbf{T}}^{t_0-\delta t}_{t_0+t}
\circ {\mathbf{T}}^{t_0+t}_{t_0}.
$$
 
\begin{thm}\label{Rt}
There exists $\delta\in ]0,1[$, which depends only on $m$ and $M$
such that the operators ${\mathbf{R}}^t,\check {\mathbf{R}}^t$
have the following properties:
\begin{itemize}
\item For each $t_0\in \Rm$ and $t\in ]0,1[$, 
the finite valued functions in the images of ${\mathbf{R}}^t$ and
 $\check {\mathbf{R}}^t$ are uniformly
$C^{1,1}$. 
\item For each semi-concave function $u$, there exists
$T>0$ such that ${\mathbf{R}}^tu\leq u$ and $\check {\mathbf{R}}^tu \leq u$
for each  $t_0\in \Rm$ and $t\in ]0,T[$.
\item For each semi-convex function $u$, there exists
$T>0$ such that ${\mathbf{R}}^tu\geq u$ and $\check {\mathbf{R}}^tu \geq u$
for each $t_0\in \Rm$ and
 $t\in ]0,T[$.
\item For each $C^{1,1}$ function $u$, 
 there exists
$T>0$ such that ${\mathbf{R}}^tu=u$ and $\check {\mathbf{R}}^tu = u$
for each $t_0\in \Rm$ and
 $t\in ]0,T[$.
\end{itemize}
\end{thm}

\proof
The finite valued functions in the image of 
${\mathbf{T}}_{t_0-t}^{t_0+\delta t}$ are $C/t$-semi-concave,
by Lemma \ref{semiconst} (we assume that  $t\in ]0,1[$).
Then, by Proposition \ref{semiconvex},
the  finite valued functions in the image of 
$\check {\mathbf{T}}_{t_0+\delta t}^{t_0}
\circ {\mathbf{T}}_{t_0-t}^{t_0+\delta t}
$
are $(2C/t)$-semi-concave provided 
$$
\delta t\leq \frac{C}{tM(3+2C/t)^2}=\frac{Ct}{M(3t+2C)^2},
$$
which holds it $\delta \leq C/(M(3+2C))$. For such a $\delta$,
the finite valued functions in the image of ${\mathbf{R}}^t$ are uniformly semi-concave.
They are also uniformly semi-convex, hence uniformly $C^{1,1}$.
The proof is similar for $\check {\mathbf{R}}$.
Let us now write 
$$
{\mathbf{R}}^t:=(\check {\mathbf{T}}_{t_0+\delta t}^{t_0}
\circ 
{\mathbf{T}}_{t_0}^{t_0+\delta t})
\circ ({\mathbf{T}}_{t_0-t}^{t_0}
\circ \check {\mathbf{T}}_{t_0}^{t_0-t}),
$$
which implies,  using the monotony, that 
 ${\mathbf{R}}^tu\geq \check {\mathbf{T}}_{t_0+\delta t}^{t_0}
\circ 
{\mathbf{T}}_{t_0}^{t_0+\delta t}u
$
and
${\mathbf{R}}^tu\leq{\mathbf{T}}_{t_0-t}^{t_0}
\circ \check {\mathbf{T}}_{t_0}^{t_0-t}u.
$
By Addendum \ref{ad} we conclude that ${\mathbf{R}}^tu\geq u$ for small $t$ when
$u$ is semi-convex. All the statements of the second and third point
follow by similar considerations.
The last point follows from the second and third one.

\qed

\section{Sub-solutions of the stationary Hamilton-Jacobi equation.}
We assume from now on that the Hamiltonian does not explicitly depend on time. Then, in addition to (\ref{HJ}), we can consider
the stationary Hamilton-Jacobi equation 
\begin{equation}\label{HJa}\tag{HJ$a$}
H(q,du(q))=a,
\end{equation}
for each real parameter $a$. This stationary equation is the main character of Fathi's joined lecture.
Formally, a function $u(q)$ solves (\ref{HJa}) if and only if 
the function $(t,q)\lmto u(q)-at$ solves (\ref{HJ}).
It is not hard to check that this also holds in the sense of viscosity solutions:
The function $u(q)$ is a viscosity solution 
of (\ref{HJa}) if and only if the function 
$(t,q)\lmto u(q)-at$ is a viscosity solution of  (\ref{HJ}).
Let us explicit for later references:
\begin{hyp}\label{h4}
We say that  $H$ is autonomous if it does not depend on the time variable.
\end{hyp}

In this autonomous context, we have 
${\mathbf{T}}_{\tau}^{\tau+t}={\mathbf{T}}^{t}_0$. We will denote by
 ${\mathbf{T}}^t$ this operator. 
The Markov property turns to the equality 
${\mathbf{T}}^t\circ {\mathbf{T}}^s={\mathbf{T}}^{t+s}$.
In other words, the Lax Oleinik operators form 
a semi-group, the famous Lax-Oleinik semi-group.
Another important specificity of the autonomous context
is that the Hamiltonian $H$ is constant along Hamiltonian
orbits, as can be checked by an easy computation.

\begin{prop}\label{sub}
Given a Hamiltonian $H$ satisfying Hypotheses \ref{h1},\ref{h2},\ref{h3},\ref{h4},
the following properties are equivalent for a function $u$:
\begin{enumerate}
\item The function $u$ is Lipschitz and it solves the 
inequation $H(q,du(q))\leq a$ almost everywhere.
\item The inequality $u(q_1)-u(q_0)\leq A^t(q_0,q_1)+at$ 
holds for each $q_0\in \Rm^d,q_1\in \Rm^d,t>0$.
\item The inequality $u\leq {\mathbf{T}}^tu+t a$ holds for each $t\geq 0$.
\item
The function $u$ is a viscosity sub-solution of the Hamilton-Jacobi equation $H(q,du(q))=a$.
\item The function $u$ is Lipschitz and  the 
inequation $H(q,du(q))\leq a$ holds at each point of differentiability $q$ of $u$ (by Rademacher Theorem,
the set of points of differentiability has full measure).
\end{enumerate}
The function $u$ is called a sub-solution at level $a$,
or a sub-solution of (\ref{HJa}), if it satisfies these properties.
\end{prop}

\proof
It is tautological that $5\Rightarrow 1$ and easy
that  $2\Leftrightarrow 3$. Let us prove that $1\Rightarrow 2$, following 
Fathi.
If $1$ holds, then there exists a set $M\subset \Rm^d$ of full measure composed of points
of differentiability $q$ of  $u$ such that $H(q,du(q))\leq a$.
We first assume that  $t< \sigma$ and prove 2 (recall that  $A^t=S^t$).
Let us consider the map
$$(q_0,q_1,\tau )\lmto (q(\tau),q_1,\tau ),
$$
where $q(\tau)$ is the value at time $\tau$ of the unique orbit
$(q(s),p(s))$ which satisfies $q(0)=q_0$ and $q(t)=q_1$.
This map 
is a diffeomorphism of $\Rm^d\times \Rm^d\times ]0,t[$, the inverse
diffeomorphism being
$$
(\theta,q_1,\tau)\lmto (q(0),q_1,\tau),
$$
where $(q(s),p(s))$ is the unique orbit such that $q(\tau)=\theta$
and $q(t)=q_1$.
As a consequence, for almost each pair $(q_0,q_1)$, the function 
$u$ is differentiable at the point 
$q(s)$ for almost every $s\in ]0,t[$.
If $(q_0,q_1)$ is such a pair, we have, using the convexity of $H$ in $p$,
\begin{align*}
u(q_1)-u(q_0)&=u(q(t))-u(q(0))=\int_0^t du_{q(s)}\cdot \dot q(s)ds=\int_0^t du_{q(s)}\cdot \partial_pH(q(s),p(s))ds\\
&\leq \int_0^t H(q(s),du_{q(s)})+\partial_pH(q(s),p(s))\cdot p(s)
-H(q(s),p(s))ds\\
&\leq at+S^t(q(0),q(t))=at+A^t(q_0,q_1).
\end{align*}
We have proved the desired inequality for almost every pair $(q_0,q_1)$, hence on a dense subset of pairs.
Since both sides of the inequality are continuous, we deduce
that the inequality holds for all pairs $(q_0,q_1)$, provided
$t<\sigma$.
In order to deduce the inequality when $t\geq \sigma$, we write,
for $n$ large enough,
\begin{align*}
A^t(q_0,q_1)+at&=\min_{\theta_1,\ldots,\theta_{n-1}}
\big(
S^{t/n}(q_0,\theta_1)+at/n+\cdots+S^{t/n}(q_{n-1},q_1)+at/n
\big)\\
&\geq \min_{\theta_1,\ldots,\theta_{n-1}}
\big(u(\theta_1)-u(q_0)+\cdots+u(q_1)-u(\theta_{n-1})
\big)=u(q_1)-u(q_0).
\end{align*}
Let us now prove that $3\Rightarrow 4$.
Let $u$ be a function satisfying 3. This function then satisfies
2, hence it is Lipschitz.
We consider a $C^2$ function $v(q)$ which touches $u$ 
from above at some point $\theta$, which means that 
$v-u$ has a global minimum  at $\theta$. Since the function $u$
is Lipschitz, 
we can modify $v$  at infinity and assume that 
it has bounded second differential.
Then, there exists a $C^2$ solution $V(t,q)$ of (\ref{HJ})
defined on  $]-T,T[\times \Rm^d$ with $T>0$, and such that 
$V(0,q)=v(q)$.
For $t\geq 0$, we have $V_t={\mathbf{T}}^tv$, by Proposition \ref{formel}.
Since $v\geq u$, we obtain that  
$$
V(t,q)={\mathbf{T}}^tv(q)\geq {\mathbf{T}}^tu(q)\geq u(q)-at
$$
for $t\in ]0,T[$, 
hence $\partial_tV(0,\theta)\geq -a$ (recall that $\theta$
 is the point of contact between $u$ and $v$).
Since we know that $V$ solves (\ref{HJ}), we conclude that 
$$
H(\theta, \partial_qV(0,\theta))=H(\theta, dv(\theta))
\leq a.
$$
The proof that $4 \Rightarrow 5$ is very classical and can be
found in Fathi's lecture, but we recall it here for completeness.
If $q$ is a point of differentiability of $u$, then $du(q)$
is a super-differential (but not necessarily a proximal super-differential) of $u$ at $q$, hence $H(q,du(q))\leq a$.
We will now prove that the function $u$ is locally Lipschitz. The estimate $H(q,du(q))\leq a$ , which holds at each point of differentiability of $u$, 
then implies that it
is globally Lipschitz in view of Hypothesis \ref{h3}.

Let $B(Q,1)$ be a closed ball, of radius one.
Let us set $r=\max_{\theta\in B(Q,2),q\in B(Q,1)}(u(\theta)-u(q))$.
Let $k$ be a positive number greater that $r$ and such that 
$|p|\geq k\Rightarrow H(q,p)>a$ for each $q$.
 Such a $k$ exists by Hypothesis \ref{h3}.
Given $q$ in $B(Q,1)$, the function 
$$
\theta\lmto k|\theta-q|-u(\theta)
$$
 has then a local minimum in the interior of the ball
$B(Q,2)$.
If this minimum is reached at a point $q_1$ different from $q$,
then the function $v(\theta):= k|\theta-q|$ is smooth at $q_1$, and, since $u$
is a viscosity sub-solution, we have $H(q_1,dv(q_1))\leq a$,
which is in contradiction with the fact that $|dv(q_1)|=k$.
Hence the minimum must be reached at $q$, which implies that 
$k|\theta-q|-u(\theta)\geq -u(q)$ or equivalently 
that 
$$
u(\theta)-u(q)\leq k|\theta-q|
$$
for each $\theta\in B(Q,2)$ and all $q\in B(Q,1)$.
We conclude that $u$ is $k$-Lipschitz on $B(Q,1)$.
\qed

\begin{cor}
If $u$ is a sub-solution of (\ref{HJa}), then, for each $t\geq 0$,
${\mathbf{T}}^tu$ is a sub-solution of (\ref{HJa}), and so is $\check {\mathbf{T}}^tu$.
\end{cor}

\proof
The function $u$ is a sub-solution if and only if
${\mathbf{T}}^su+as\geq u$ for each $t\geq 0$.
Applying ${\mathbf{T}}^t$, we obtain 
${\mathbf{T}}^t{\mathbf{T}}^su+as={\mathbf{T}}^s{\mathbf{T}}^tu+as\geq {\mathbf{T}}^tu$.
Since this inequality holds for each $s\geq 0$, we conclude
that ${\mathbf{T}}^tu$ is a sub-solution. 
\qed

\begin{cor}
If the function $u$ is Lipschitz, and if the Hamiltonian is autonomous, then the functions ${\mathbf{T}}^tu,t\geq 0$
are equi-Lipschitz.
\end{cor}
\proof
If the function $u$ is $k$-Lipschitz, then $du(q)\leq k$ almost  everywhere, hence $u$ is a sub-solution to (\ref{HJa})
for some $a$ (one can take $a=\sup_{|p|\leq k}H(q,p)$).
As a consequence, the functions ${\mathbf{T}}^tu,t\geq 0$
are all sub-solutions to (\ref{HJa}), hence they are $K$-Lipschitz,
with $K=\sup \{|p|, H(q,p)\leq a\}$.
\qed

\section{Weak KAM solutions and invariant sets.}\label{WKAM}
We derive here the first dynamical consequences from the theory.
\begin{defn}
The function $u$ is called a Weak KAM solution at level $a$ if 
${\mathbf{T}}^tu+ta =u$
for each $t\geq 0$. 
Weak KAM solutions at level $a$ are viscosity solutions of 
(\ref{HJa}).
We say that the function $u$ is a Weak KAM Solution if it is 
a Weak KAM solution at some level $a$.
\end{defn}

If $u$ is a weak KAM solution, then it is semi-concave
(with a semi-concavity constant which depends only on 
$M$ and $m$).
By Theorem \ref{geom}, for $t>0$, we have the inclusion 
$$
\varphi^{-t}\big(\bar \Gamma_u\big)\subset \Gamma_u
$$
and this set is a Lipschitz graph.
The set 
$$
\mI^*(u):=\cap_{n\in \Nm} \varphi^{-n}\big(\bar \Gamma_u\big)
$$
is a closed invariant set contained in a Lipschitz graph.
It would be a very nice result to have obtained a 
distinguished closed
invariant subsets of our Hamiltonian system contained in 
a Lipschitz graph. Unfortunately, at this point, we can't
prove (because it is not necessarily true) that 
the set $\mI^*(u)$ is not empty.
In order  to  obtain interesting dynamical
consequences from this theory, we need an 
additional assumption.

\begin{hyp}\label{h5}
We say that the Hamiltonian $H$ is periodic if $H(q+w,p)=H(q,p)$
for each $w\in \Zm^d$, $q\in \Rm^d$ and $p\in \Rm^{d*}$.
\end{hyp}

Under this hypothesis, we should see the Hamiltonian system
as defined on the  phase space
$\Tm^d\times \Rm^{d*}$, with $\Tm^d=\Rm^d/\Zm^d$.
Indeed, the flow $\varphi^t$ commutes with the translations
$(q,p)\lmto(q+w,p)$, $w\in \Zm^d$.
The compactness of this new configuration space has 
remarkable consequences, summed up in the following Theorem.
We assume in the rest of this section that the Hamiltonian 
$H$ satisfies Hypotheses \ref{h1}, \ref{h2}, \ref{h3}, \ref{h4}, \ref{h5}.

\begin{thm}
If the Hamiltonian is autonomous and periodic, then there
exists a periodic Weak KAM solution.
The corresponding set $\mI^*(u)$ is a non-empty closed 
invariant set which is contained in a Lipschitz
graph and which is invariant under the translations
$(q,p)\lmto(q+w,p)$, $w\in \Zm^d$.  
\end{thm}
This last property on the invariance under translations means that 
$\mI^*(u)$ naturally gives rise to an invariant space on the 
quotient phase space $\Tm^d\times \Rm^{d*}$.

\proof
Let us first prove the second part of the Theorem.
If $u$ is a periodic Weak KAM solution, then the 
set $\bar \Gamma_u$ is 
contained in $\{|p|\leq C\}$ for some constant $C$, and it is 
invariant under the integral translations,
hence it descends to a compact subset of $\Tm^d\times \Rm^{d*}$,
that we still denote by $\bar \Gamma_u$.
Then the sets $\varphi^{-n}(\bar \Gamma_u)$ form a decreasing
sequence of non-empty  compact sets, hence their intersection is a 
non-empty compact set.

Let us now prove that there exists a periodic Weak KAM solution.
We follow the proof of \cite{JAMS}, which is slightly different
from the original proof of Fathi.
Observe first that the functions $A^t(q_0,q_1)$
are periodic in the sense that $A^t (q_0+w,q_1+w)=A^t(q_0,q_1)$
for each $w\in \Zm^d$.
This implies that ${\mathbf{T}}^tu$ is periodic when $u$ is periodic.
Considering the Cauchy problem for (\ref{HJ})
with initial condition equal to zero,
we define $v(t,q):={\mathbf{T}}^t0(q)$.
The quantities  $a^+(t)=\max_q v_t(q)$
and $a^-(t)=\min_q v_t(q)$ will be useful.
Since the functions $v_t, t\geq 0$ are equi-Lipschitz, 
there exists a constant $K$ such that $a^+(t)-a^-(t)\leq K$
for all $t\geq 0$.
We have
$$
a^+(t+s)=\max {\mathbf{T}}^{t+s}0=\max {\mathbf{T}}^t({\mathbf{T}}^s0)\leq 
{\mathbf{T}}^t(a^+(s))=a^+(s)+{\mathbf{T}}^t(0)\leq a^+(s)+a^+(t),
$$
and similarly 
$$
a^-(t+s)\geq a^-(t)+a^-(s).
$$
By standard results on sub-additive functions, we conclude that 
$a^+(t)/t$ and $a^-(t)/t$ converge respectively to 
$\inf_{t\geq 0} a^+(t)/t$ and $\sup_{t\geq 0}a^-(t)/t$.
Since $a^+-a^-$ is bounded, these two limits have the same value, 
let us call it $-a$.
We have 
$$
K-ta\geq a^-(t)+K\geq a^+(t)\geq -ta\geq a^-(t)\geq a^+(t)-K\geq K-at
$$
for all $t\geq 0$, hence
$$
K\geq v(t,q)+ta\geq -K.
$$
We can now define 
$$
u(q):=\liminf _{t\lto \infty} (v(t,q)+ta).
$$
We claim that $u$ is a Weak KAM solution at level $a$.
Since the functions $v_t+ta$ are equi-Lipschitz and equi-bounded,
the function $u$ is well-defined and Lipschitz.
We have to prove that ${\mathbf{T}}^tu+ta=u$ for all $t\geq 0$.

We have 
$$
v(t+s,q_1)+(t+s)a\leq v(s,q_0)+sa+A^t(q_0,q_1)+ta
$$
for each $q_0, q_1$ and $t\geq0, s\geq 0$.
Taking the $\liminf$ in $s$  yields
$$
u(q_1)\leq u(q_0)+A^t(q_0,q_1)+ta.
$$
We have proved that $u$ is a sub-solution to (\ref{HJa}).

Conversely, we have to prove that ${\mathbf{T}}^tu+ta\geq u$.
Let us pick a point $q$ and consider a sequence $t_n$ such 
that $v(t_n,q)+t_na\lto u(q)$. Fixing $t>0$, 
we consider a sequence  $q_n$  in $\Rm^d$ such that 
$$v(t_n,q)+t_n a=v(t_n-t,q_n)+(t_n-t)a+A^t(q_n,q)+ta.
$$
This equality implies that the sequence $q_n$ is bounded,
and we  assume by taking a subsequence that it has a limit
$q'$. We can also assume that the sequence $v(t_n-t,q')+(t_n-t)a$
has a limit, that we denote by $l$. Note that $l\geq u(q')$.
Since the functions $v_t$ are equi-Lipschitz, we have
$v(t_n-t,q_n)+(t_n-t)a\lto l$
hence, taking the limit in the equality above,
$$
u(q)=l+A^t(q',q)+at\geq u(q')+A^t(q',q)+at\geq {\mathbf{T}}^tu(q)+at.
$$ 
We have proved that $u$ is a periodic Weak KAM solution at level $a$.
\qed

The periodic Weak KAM solutions at level $a$
are the periodic viscosity solutions of (\ref{HJa}),
as is proved in Fathi's joined lecture.
The existence of periodic viscosity solutions was first
obtained by Lions, Papanicolaou and Varadhan  in a famous unpublished preprint, \cite{LPV}.
The most important aspect of Fathi's weak KAM theorem that we 
just exposed is that these viscosity solutions have a dynamical relevance and give rise to invariant sets.

Let us comment a bit further in that direction, and explain the 
name ''Weak KAM``.
Consider a periodic Lipschitz  function $u$, and the associated
set $\Gamma_u$, that we consider here as a subspace of 
$\Tm^d\times \Rm^{d*}$.

Assume first that $u$ is $C^2$, so that $\Gamma_u$ is a $C^1$
graph. This graph is invariant if and only if there 
exists $a$ such that $u$ solves (\ref{HJa}).
This follows from Section \ref{L1}: If $u$ solves (\ref{HJa}),
then the function $U(t,q)=u(q)-at$ solves $\ref{HJ}$,
hence 
$$
\varphi^t(\Gamma_{u})=\Gamma_{U_t}=\Gamma_{u}.
$$
Conversely, if $\Gamma_u$ is invariant, then 
$\Gamma _{{\mathbf{T}}^tu}=\varphi^t(\Gamma_u)=\Gamma_u,$ by Corollay \ref{flow},  hence 
${\mathbf{T}}^tu$
is equal to $u$ up to an additive constant $a(t)$.
Since ${\mathbf{T}}^t$ is a semi-group, it is easy to deduce that $a(t)=at$
for some $a\in \Rm$. As a consequence,  $u$ is a $C^2$ Weak KAM solution, hence  a classical solution of (\ref{HJa}).

The classical KAM theorem gives the existence, in certain very
specific settings, of some invariant $C^1$ graphs of the form
$\Gamma_u$. From what we just explained, it can be interpreted as
giving the existence of $C^2$ solutions of (\ref{HJa}),
although this point of view is not the right one to obtain its proof.
It is  natural to expect that the Hamilton-Jacobi equation
could be used to produce invariant sets in more general situations.
Since we do not know any direct method to prove the existence of 
$C^2$ solutions of (\ref{HJa}), we should deal with 
some kind of weak solutions.
However, if $u$ is just  a Lipschitz solution almost everywhere,
we can't say much about the dynamical properties 
of $\Gamma_u$.
It is remarkable that the inclusion $\varphi^t( \Gamma_u)\supset \bar \Gamma_u$
 holds for viscosity solutions (or, equivalently
Weak KAM solutions) in the convex case.
This is the starting point of 
Fathi's construction of the  invariant set $\mI^*(u)$
that  we  exposed in the present section.

\section{Regular sub-solutions and the Aubry set.}\label{regular}

We abandon for a moment the hypothesis \ref{h5},
and consider a Hamiltonian satisfying Hypotheses \ref{h1}, \ref{h2}, \ref{h3}, \ref{h4}.
We describe a new construction of invariant sets 
based on the study of regular sub-solutions, and define the Aubry set.
We mostly follow \cite{ENS} in this section.
The following result is at the base of our constructions,
see \cite{ENS,Ar:11,FS:04}.

\begin{thm}\label{C11}
If (\ref{HJa}) admits a sub-solution, then it admits a
$C^{1,1}$ sub-solution.
Moreover, the set of $C^{1,1}$ sub-solutions is dense in the 
set of all sub-solutions for the uniform topology.
\end{thm}

\proof
Let $u$ be a sub-solution at level $a$.
We use the operator
 $\mathbf{R}^t=\check {\mathbf{T}}^{\delta t} 
\circ  {\mathbf{T}}^{(\delta+1) t}
\circ \check {\mathbf{T}}^{ t}$ 
of Theorem \ref{Rt} to regularize $u$.
Since the operators $\mathbf{T}^t$ and $\check {\mathbf{T}}^t$
preserve sub-solutions, so does $\mathbf{R}^t$.
We claim that 
$$
 u-(C+a)(1+\delta)t \leq  {\mathbf{R}}^t u\leq u+(C+a)(1+\delta)t
$$
with a constant $C$ which depends only on $m$ and $M$.
This  implies that the function ${\mathbf{R}}^t u$
is finite valued.
If the parameter $\delta$ has been chosen small enough,
then, by Theorem \ref{Rt},  the functions $\mathbf{R}^t$ are $C^{1,1}$ sub-solutions,
which converge uniformly to $u$ as $t\lto 0$.
The bound on ${\mathbf{R}}^t u$ claimed above follows from
the following ones in view of Property \ref{axiom}, 
$$
v-sa \leq {\mathbf{T}}^sv \leq v+Cs
,\quad 
v-Cs \leq \check {\mathbf{T}}^s v\leq v+sa
$$
which hold for each $s\geq 0$ and each sub-solution $v$ at level $a$. 
The first one can be seen by writing 
$$
u(q)-as\leq {\mathbf{T}}^s u(q)\leq u(q)+A^s(q,q)\leq u(q)+Cs.
$$
This ends the proof of Theorem \ref{C11}.
Observe that we could have used the simpler operator
$\breve {\mathbf{T}}^{\delta t}\circ  {\mathbf{T}}^t$,
as was done in \cite{ENS}, but the operator $\mathbf{R}^t$
deserves attention for some nicer properties. 
\qed

\begin{defn}
The critical value of $H$ is the real number $\alpha$
(or $\alpha(H)$) defined as  the infimum of all real numbers 
$a$ such that (\ref{HJa}) has a sub-solution.
The sub-solutions of (HJ$\alpha$) are called critical sub-solutions. 
\end{defn}

\begin{lem}
We have the estimate $-M\leq \alpha\leq M$.
\end{lem}
\proof
The function $u=0$ is a sub-solution at level $M$, hence 
$\alpha\leq M$.
Conversely, since $H\geq -M$ there exists no sub-solution
at level $a$ when $a<-M$.
\qed

\begin{prop}
There exists a $C^{1,1}$ sub-solution of (HJ$\alpha$).
\end{prop}
\proof
Let $a_n$ be a sequence decreasing to $\alpha$.
Since $a_n>\alpha$, the Hamilton-Jacobi equation at level $a_n$
has a sub-solution $u_n$.
The sequence $u_n$ is equi-Lipschitz, 
and we can assume by adding constants that it is also equi-bounded.
Taking a subsequence,
we can also assume  that it converges locally uniformly to a limit $u$.
Taking the limit $n\lto \infty$ in the inequalities 
$u_n(q_1)-u_n(q_0)\leq A^t(q_0,q_1)+ta_n$
gives $u(q_1)-u(q_0)\leq A^t(q_0,q_1)+ta_n$.
This holds for all $q_0,q_1$ and $t>0$, hence $u$ is a sub-solution
at level $\alpha$, or in other words a critical sub-solution.
Since there exists a critical sub-solution, Theorem \ref{C11}
implies that there exists a $C^{1,1}$ critical sub-solution.
\qed

\begin{defn}
The projected Aubry set is the set $\mA\subset \Rm^d$ of points 
$q$ such that the equality $H(q,du(q))=\alpha$ holds for all
$C^1$ critical sub-solutions $u$.
\end{defn}
We point out that $\mA$ might be empty without additional hypotheses.

\begin{lem}\label{31}
If $q\in \mA$, then all $C^1$ critical sub-solutions $u$
have the same differential at $q$. In other words, the restriction
$\Gamma_{u|\mA}$ does not depend on the $C^1$ critical sub-solution 
$u$.
\end{lem}

\proof
Let $u$ and $v$ be two critical sub-solutions, and $q$ 
a point in $\mA$. We have to prove that $du(q)=dv(q)$.
Assume, by contradiction, that this equality does not hold
and consider the sub-solution $w=(u+v)/2$. 
Since $H(q,du(q))=H(q,dv(q))=\alpha$, the strict convexity 
of $H(q,.)$ implies that  $H(q,dw(q))<\alpha$, which contradicts the definition of $\mA$.
\qed

\begin{lem}
There exists a $C^{1,1}$ sub-solution $u_0$ 
which satisfies the strict inequality $H(q,du_0(q))<\alpha$ for all $q$ in the complement
of $\mA$.
\end{lem}
\proof
The set of $C^1$ functions is separable for the topology of uniform  $C^1$ convergence on compact sets.
This topology can be defined for example by the distance 
$$
d(u,v)=\sum_n \frac{\sup_{|q|\leq n} \arctan(|u(q)|+|du(q)|)}{2^n}.
$$
Since a subset of a separable space is separable,
there exists a sequence $u_n$ of $C^1$ critical  sub-solutions
which is dense for this topology in the set of all $C^1$
critical 
sub-solutions.
Let us set  
$$
a_n=\frac{a_0}{2^n \sup_{k\leq n, |q|\leq n} (1+|u_k(q)|+|du_k(q)|))}
$$
and  choose $a_0$ such that $\sum_{n\geq 1} a_n=1$.
The sum $\sum_{n\geq 1} a_nu_n$ converges uniformly
with its differentials on each compact sets to a $C^1$ limit $v_0$.
The function $v_0$ is a critical sub-solution, and we claim that 
 $H(q,dv_0(q))=\alpha$  if and only if $q$ belongs to $\mA$.
Indeed, this equality holds only if all the inequalities
$H(q,du_n(q))\leq \alpha$ are equalities, which, in view of the 
density of the sequence $u_n$, implies that $H(q,du(q))=\alpha $
for all $C^1$ sub-solutions $u$. By definition, this implies that
$q$ belongs to $\mA$.
We have constructed a $C^1$ sub-solution $v_0$ such that 
$$
H(q,dv_0(q))<\alpha
$$
outside of $\mA$. We have to prove the existence of a $C^{1,1}$
critical sub-solution with the same property.
We consider a smooth function $V(q)$ which is bounded in $C^2$,
which is positive outside of $\mA$, and
such that 
$$
0\leq V(q)\leq \alpha-H(q,dv_0(q))
$$
for all $q\in \Rm^n$.
The modified Hamiltonian $\tilde H(q,p)=H(q,p)+V(q)$
satisfies all our hypotheses.
Since $\tilde H\geq H$, the corresponding critical value
$\tilde \alpha$ satisfies $\tilde \alpha\geq \alpha$.
Since $v_0$ is a sub-solution of the inequation
$$
\tilde H(q,dv_0(q))\leq \alpha,
$$
we can apply Theorem \ref{C11} to $\tilde H$ at level $\alpha$,
and obtain the existence of a $C^{1,1}$  sub-solution $u_0$ to the same inequation.
The inequality 
$$
H(q,du_0(q))\leq \alpha-V(q)
$$
implies that $u_0$ is a critical sub-solution for $H$
which is strict on the set $\{V>0\}$ which, from 
our construction of $V$, is the complement of $\mA$.
\qed
\begin{defn}
 The Aubry set  $\mA^*$ is defined as:
$$
\mA^*=\cap _u\Gamma_{u|\mA}=\cap _u \Gamma_u,
$$
where the intersections are taken on the set of $C^1$ critical sub-solutions.
\end{defn}

In view of Lemma \ref{31} we have $\mA^*=\Gamma_{u|\mA}$ 
for each $C^1$ sub-solution $u$, hence $\pi(\mA^*)=\mA$, 
where $\pi:\Rm^d\times \Rm^{d*}\lto \Rm^d$
is the projection on the first factor. 
To check the second inequality, it is sufficient to prove that 
$\cap_u\Gamma_u\subset \mA^*$.
Let $u_0$ be a $C^1$ critical  sub-solution such that $H(q,du_0(q))<\alpha$ 
outside of $\mA$.
Given a point $(q_0,p_0)$ in $\Gamma_{u_0}-\mA^*$, we can slightly
perturb the critical sub-solution $u_0$ around $q_0$ to a 
critical sub-solution $u_1$
such that $du_1(q_0)\neq du_0(q_0)$ (we use  the 
strict inequality $H(q,du_0(q))<\alpha$). 
The point $(q_0,p_0)$ does not belong to $\Gamma_{u_1}$, hence it
does not belong to $\cap_u\Gamma_u$, which ends our proof.

The set $\mA^*$ is contained in the Lipschitz graph $\Gamma_{u_0}$
for  each $C^{1,1}$ sub-solution $u_0$.
As in Section \ref{WKAM}, we have obtained an invariant
set contained in a Lipschitz graph, but which may be empty
in general:

\begin{prop}
The Aubry set is a closed invariant set. 
\end{prop}

\proof
Let $u_0$ be a $C^{1,1}$ critical solution such that 
$H(q,du_0(q))<\alpha$ outside of $\mA$.
By Proposition \ref{semiconvex}, there exists $T>0$ such that ${\mathbf{T}}^tu_0$ is still $C^{1,1}$
for $t\in [-T,T]$.
Given $(q,p)\in \mA^*$, we conclude that, for $t\in [0,T]$,
we have  $p=d({\mathbf{T}}^tu_0)(q)$.
Setting $\theta=Q^{-t}(q,p)$,
the addendum to Proposition \ref{semiconvex} implies that
${\mathbf{T}}^tu_0(q)=u_0(\theta)+S^t(\theta,q)$,
and that 
$$
\varphi^t(\theta,du_0(\theta))=
(q,p).
$$
Since the flow preserves the Hamiltonian,
we get that
$
H(\theta,du_0(\theta))=\alpha,
$
hence the point $\theta$ belongs to $\mA$, and then 
$$
\varphi^{-t}(q,p) =\big (\theta,du_0(\theta)\big)\in \mA^*.
$$
We have proved that $\varphi^{-t}(\mA^*)\subset \mA^*$
for $t\in [0,T]$.
We can prove in a similar way, using the $C^{1,1}$ sub-solution 
$\check {\mathbf{T}}^tu_0$ instead of ${\mathbf{T}}^tu_0$, that
$\varphi^{t}(\mA^*)\subset \mA^*$
for $t\in [0,T]$, and hence that 
$$
\varphi^{t}(\mA^*)=\mA^*
$$
for each $t\in [-T,T]$, which clearly implies that this equality 
holds for all $t$.
We have proved the invariance of $\mA^*$.
\qed

\begin{prop}\label{aubrylax}
The equality
$$
\check {\mathbf{T}}^tu(q)-t\alpha=u(q)={\mathbf{T}}^tu(q)+t\alpha
$$
holds for each critical sub-solution $u$, each $t\geq 0$ and each 
$q\in \mA$.
The inclusion $\mA^*\subset \Gamma_u$ holds
for each critical sub-solution,
hence the inclusion $\mA^*\subset \mI^*(u)$ holds for each
 weak KAM solution at level $\alpha$.
\end{prop}

\proof
Let $(q(s),p(s))$ be a trajectory contained in $\mA^*$, 
 and $t\geq 0$ be given.
 For each   $C^1$ 
critical sub-solution $u$,
we have $p(s)=du_{q(s)}$, and
\begin{align*}
u(q(t))-u(q(0))=&
\int_0^t du_{q(s)}\dot q(s) ds
=
t\alpha+\int_0^t du_{q(s)}\dot q(s)-H(q,du_{q(s)})ds\\
\geq
&A^{t}(q(0),q(t))+t\alpha.
\end{align*}
Since $u$ is a critical sub-solution,  
the second point in  Proposition \ref{sub}
 implies that  the last inequality must be an equality, hence
$$
u(q(t))-u(q(s))=A^{t-s}(q(s),q(t))+(t-s)\alpha
$$
for each $t\geq s$.
In the terminology of Fathi, we have proved that the curve $q(s)$
is calibrated by the sub-solution $u$.
We can now write 
$$
u(q(t))\leq {\mathbf{T}}^tu(q(t))+t\alpha 
\leq u(q(0))+A^t(q(0),q(t))+t\alpha=u(q(t)).
$$
This implies that ${\mathbf{T}}^tu+t\alpha=u$ on $\mA$, and, similarly, 
$\check {\mathbf{T}}^tu-t\alpha =u$
on $\mA$.
Let us now fix  $t\in ]0,\sigma[$.
Given an orbit $(q(s),p(s))$ in $\mA^*$, we have 
$$
u(q(0))\leq u(\theta)+S^t(\theta,q(0))+t\alpha
$$
for each sub-solution $u$ and 
 each $\theta$, with equality at $\theta=q(-t)$.
This implies that $\partial_1S(q(-t),q(0))$
is a super-differential of $u$ at $q(0)$.
This holds in particular for $C^1$ sub-solutions, which satisfy 
$du(q(0))=p(0)$, hence $\partial_1S(q(-t),q(0))=p(0).$
We have proved that $p(0)$ is a super-differential of 
$u$ at $q(0)$. Similarly, using the inequality
$$
u(q(0))\geq u(\theta)-S^t(q(0),\theta)-t\alpha,
$$
with equality at $\theta=q(t)$,
we conclude that $p(0)$ is a sub-differential of $u$ at $q(0)$.
This implies that $u$ is differentiable at $q(0)$, and that
$du(q(0))=p(0)$. As a consequence, $\mA^*\subset \Gamma_u$
for each sub-solution $u$.
\qed

In the course of the above proof, we have established the following lemma, which will be needed later:

\begin{lem}\label{Tu=u}
Let $u$ be a sub-solution at level $a$, and let $(q(s),p(s))$ be a Hamiltonian
trajectory contained in $\Gamma_u\cap\{H=a\}$ (note that
this set is not necessarily invariant in general),
then, the equality
$\check {\mathbf{T}}^tu(q(s))-ta=u(q(s))={\mathbf{T}}^tu(q(s))+ta$ holds, for each $t\geq 0$
and  each $s\in \Rm$.
\end{lem}

\section{The Ma\~n\'e Potential.}
 In this section, we work with a Hamiltonian satisfying Hypotheses 
\ref{h1}, \ref{h2}, \ref{h3}, \ref{h4}.
The Ma\~n\'e Potential at level $a$ is the function 
$$
\Phi^a(q_0,q_1):= \inf_{t>0}\big( A^t(q_0,q_1)+at
\big).
$$
This function was first introduced by Ricardo Ma\~n\'e, see \cite{Mane:97}.
We leave as an easy exercise for the reader to prove the triangle inequality
$$
\Phi^a(q_0,q_1)\leq \Phi^a(q_0,\theta)+\Phi^a(\theta,q_1).
$$
In view of Proposition \ref{sub}, each sub-solution $u$
at level $a$ satisfies 
$$
u(q_1)-u(q_0)\leq \Phi^a(q_0,q_1)
$$
for each $q_0$ and $q_1$.
We conclude that $\Phi^a$ is finite if there exists a sub-solution
at level $a$, which holds if and only if $a\geq \alpha$.
Conversely, If the function $\Phi^a$ is finite, 
then we see from the triangle inequality that the function
$q\lmto \Phi^a(q_0,q)$
is a sub-solution at level $a$, which implies that $a\geq \alpha$.
The estimates of Lemma \ref{estimeA} imply that 
$$
\Phi^a(q_0,q_1)\leq 2\sqrt{2m(M+a)}|q_1-q_0|
$$
provided $a\geq \alpha$ 
(note that $\alpha\geq -M$ ).
We have proved  that the Ma\~n\'e Potential is the function 
called the viscosity semi-distance in Fathi's lecture:
\begin{prop}
If $a\geq \alpha$, then the function 
$q\lmto \Phi^a(q_0,q)$ is the 
maximum  of all sub-solutions $u$ at level $a$ which vanish 
at $q_0$.
If $a<\alpha$, then there is no such sub-solution and
$\Phi^a$ is identically equal to $-\infty$.
\end{prop} 

This statement also implies that the Ma\~n\'e Potential
 at level $a$ only depends on the  energy level
$\{H=a\}$. More precisely, let $G$ be  another Hamiltonian
satisfying our hypotheses and such that 
$H=a\Leftrightarrow G=a$.
Then, the sets $\{H\leq a\}$ and $\{G\leq a\}$ are equal,
which implies in view of 
the first characterization of sub-solutions in Proposition \ref{sub}
that  $G$ and $H$ have the same
sub-solutions at level $a$.
As a consequence, they have the same Ma\~n\'e potential
at level $a$.
 This  is also reflected in the following 
Proposition by the fact that the involved orbits
are contained in the set  $\{H=a\}$.

\begin{prop}\label{calibree}
Given $q_0 \neq q_1$ , 
there exists $\tau \in ]0,\infty]$ and an orbit
$(q(s),p(s)):(-\tau,0]\lto \Rm^d\times \Rm^{d*}$
such that $q(0)=q_1$,
 $A_s^0(q_0,q(s))-as=\Phi^a(q(s),q_1)$,
$$
\Phi^a(q_0,q(s))+\Phi^a(q(s),q_1) =\Phi^a(q_0,q_1)
$$
and $H(q(s),p(s))= a$ for each $s\in (-\tau,0]$.
If moreover $\tau$ is finite, then $q(-\tau)=q_0$.
\end{prop}

\proof
If $q_0\neq q_1$, then either the functions
$t\lmto A^t(q_0,q_1)+at$ reaches its minimum
at some finite time $\tau>0$, or it has a 
minimizing sequence $\tau_n\lto \infty$. This follows from Lemma \ref{estimeA}.

\textbf{In the first case}, there exists an orbit 
$(q(t),p(t)):[-\tau,0]\lto \Rm^d\times \Rm^{d*}$
such that $q(-\tau)=q_0$, $q(0)=q_1$, and
$$
\int_{-\tau}^0 p\cdot \dot q -H(q,p) dt =
A^{\tau}(q_0,q_1)=\Phi^a(q_0,q_1)-\tau a.
$$
We obtain, for each $s\in [-\tau,0]$, that 
\begin{align*}
\Phi^a(q_0,q_1)-a\tau=\int_{-\tau}^0
p\cdot \dot q -H(q,p) dt 
&= \int_{-\tau}^s
p\cdot \dot q -H(q,p) dt +
\int_{s}^0
p\cdot \dot q -H(q,p) dt \\
&\geq A^{s+\tau}(q_0,q(s))+A^{-s}(q(s),q_1)\\
&\geq \Phi^a(q_0,q(s))-a(s+\tau)+\Phi^a(q(s),q_1)
+as\\
&\geq \Phi^a(q_0,q_1)-a\tau.
\end{align*}
We conclude that all these inequalities are equalities,
hence
$$
\Phi^a(q_0,q(s))+\Phi^a(q(s),q_1) =\Phi^a(q_0,q_1).
$$
We also deduce from the above chain of inequalities
that $A^{-s}(q(s),q_1)-as=\Phi^a(q(s),q_1)$,
which implies that 
the function 
$t\lmto
A^t(q(s),q_1)+at
$
is minimal for $t=-s$. Taking $s\in ]-\sigma,0[$,
we can differentiate with respect to $t$ at 
$t=-s$ and get 
$$
\partial_{t|t=-s} S^t(q(s),q_1)+a=0.
$$
Recalling the equality 
$$
\partial_tS^{-s}(q(s),q_1)+H(q_1,p(0))=0
$$
(because $p(0)=\rho_1(-s,q(s),q_1)$ in the notations 
of Section \ref{sectionconvex}), we conclude that
$H(q_1,p(0))=a$, and, since the Hamiltonian is constant
on Hamiltonian orbits, $H(q(t),p(t))=a$ for each $t$.

\textbf{In the second case}, there exists a sequence 
of orbits $(q_n(t),p_n(t))$ on $[-\tau_n,0]$
such that 
$$
\int_{-\tau_n}^0 p_n\cdot \dot q_n -H(q_n,p_n) dt+a\tau_n=
A^{\tau_n} (q_0,q_1) +a\tau_n\leq \Phi ^a(q_0,q_1)+\delta_n,
$$
where $\delta_n\lto 0$.
Let us denote $h_n:= H(q_n(s),p_n(s))$, it does not depend on $s$.
By Lemma \ref{estimeeaction} and the above inequality,
we have
$$
\frac{m}{M}\tau_n h_n-(m+M)\tau_n \leq 
\int_{-\tau_n}^0p_n\cdot \partial_pH(q_n,p_n)-H(q_n,p_n)dt
\leq \Phi ^a(q_0,q_1)+\delta_n
$$
 hence the sequence $h_n$ is bounded. 
As a consequence, the curves $p_n(s)$ are uniformly bounded,
hence so is 
$\dot q_n(s)=\partial_pH(q_n(s),p_n(s))$.
On each compact interval of time $[s,0]$,
the curves $x_n(t)=(q_n(t),p_n(t))$ are thus uniformly bounded, hence uniformly Lipschitz.
Up to taking a subsequence, we can thus assume that
the curves $x_n(t)$ converge, uniformly on compact
time intervals, to a Hamiltonian orbit 
$x(t)=(q(t),p(t)):(-\infty,0]\lto \Rm^d\times \Rm^{d*}$.
Passing at the limit in the inequality
$$
\Phi^a(q_0,q_n(s))+\Phi^a(q_n(s),q_1)\leq \Phi^a(q_0,q_1)+\delta_n,
$$
which holds for each $s\in [-\tau_n,0]$, yields
$$
\Phi^a(q_0,q(s))+\Phi^a(q(s),q_1)\leq \Phi^a(q_0,q_1),
$$
which must actually be an equality.
We prove as in the first case that $H(q_1,p(0))=a$,
thus $H(q(s),p(s))\equiv a$.
\qed

The projected Aubry set $\mA$ can be characterized in 
terms of the Ma\~n\'e potential (see also Fathi's lecture):
\begin{prop}
The following statements are equivalent for a point 
$q_0$ and a real number $a$, where we denote by $u$ the function $\Phi^a(q_0,.)$:
\begin{enumerate}
\item
$q_0\in \mA$ and $a=\alpha$.
\item
${\mathbf{T}}^tu(q_0)+ta=u(q_0)=0$ for each $t\geq 0$.
\item
The function $u$ is a Weak KAM solution at level $a$.
\item
$u$ is differentiable at $q_0$.
\end{enumerate}
\end{prop}
\proof
$1\Rightarrow 2$. 
This follows from Proposition \ref{aubrylax} since $u$ is a sub-solution 
at level $a=\alpha$.

$2\Rightarrow 3$.
Let us fix $t>0$ and $q_1$. We have to prove that 
there exists $\theta$ such that 
$u(q_1)\geq u(\theta)+A^t(\theta,q_1)+ta$
(this inequality is then an equality).
If $q_1=q_0$, the existence of this point follows from
the equality ${\mathbf{T}}^tu(q_0)+ta=u(q_0)$.

If $q_1\neq q_0$, we can apply Proposition \ref{calibree} to this pair of points.
With the notations of Proposition \ref{calibree}, 
if $\tau\geq t$, then the point $\theta=q(-t)$ fulfills
our demand. If $\tau<t$, then we set $s=t-\tau$.
We have 
$q(-\tau)=q_0$ and $A^{\tau}(q_0,q_1)+a\tau=u(q_1).$  
Since ${\mathbf{T}}^su(q_0)+sa=u(q_0)$, there exists $\theta$ such that 
$u(\theta)+A^s(\theta,q_0)+sa=u(q_0)=0$. The infimum in the definition
of  ${\mathbf{T}}^su(q_0)$ exists because $u$ is Lipschitz.
We conclude that 
$$
u(\theta)+A^t(\theta,q_1)+at\leq
 u(\theta)+A^s(\theta,q_0)+sa + A^{\tau}(q_0,q_1)+a\tau=u(q_1).
$$

$3\Rightarrow 4$.
If $u$ is a Weak KAM solution, then it has a proximal super-differential at each point. Conversely, if 
$v$ is a $C^1$ sub-solution, then $u-v$ has a minimum
at $q_0$ hence $dv(q_0)$ is a sub-differential of $u$
at $q_0$. The function $u$ both has a super-differential
and a sub-differential at $q_0$, hence it is differentiable at $q_0$.

$4\Rightarrow 1$.
If $a>\alpha$ or if $q_0$ does not belong to
$\mA$, then there exists a $C^1$ sub-solution $v$ at level
$a$ which is strict near $q_0$. 
We can then slightly perturb the function $v$
near $q_0$ and build a sub-solution $w$
such that $dw(q_0)\neq dv(q_0)$.
In view of the characterization of $u$ as the largest
sub-solution vanishing at $q_0$, we conclude
that $dv(q_0)$ as well as $dw(q_0)$ are sub-differentials
of $u$ at $q_0$, hence $u$ is not differentiable at this point.
\qed

The Ma\~n\'e potential also allows to build Weak KAM solutions
in the non periodic case by the Busemann method,
see \cite{C:01} and
 Fathi's Lecture.
Let $q_n$ be a sequence of points of $\Rm^d$ such that
$|q_n|\geq n$.
We consider the sequence of functions 
$$
u_n(q)=\Phi^{a}(q_n,q)-\Phi^a(q_n,q_0).
$$
By construction, $u_n(q_0)=0$, and it follows from
the triangle inequality that the functions $u_n$
are equi-Lipschitz. We can then assume, without loss
of generality, that the functions $u_n$ converge, 
uniformly on compact sets, to a Lipschitz limit $u(q)$.

\begin{prop}
The limit function $u(q)$ is a Weak KAM solution at level $a$.
\end{prop}
\proof
The functions $u_n$ are all sub-solutions at level $a$,
which means that $u_n(q_1)-u_n(q_0)\leq A^t(q_0,q_1)+ta$
for each $t\geq 0$, $q_0$, $q_1$.
At the limit $n\lto \infty$, we obtain that  that ${\mathbf{T}}^tu+ta\geq u$
for each $t\geq 0$.

We have to prove that ${\mathbf{T}}^tu+ta\leq u$ for all $t\geq 0$.
Let us fix a point $q$ and a time $t\geq 0$,
and consider a sequence $t_n$ such that 
$$
A^{t_n}(q_n,q)+at_n\leq \Phi^a(q_n,q)+1/n.
$$
This inequality implies that 
$$
\frac{1}{2Mt_n}|q_n-q|^2\leq 1+(M-a)t_n+2\sqrt{2m(M+a)}|q_n-q|
$$
and, since $|q_n-q|\lto \infty$, that $t_n\lto \infty$.
When $n$ is large enough, we have $t_n\geq t$ and there 
exists $\theta_n\in\Rm^d$ such that
$A^{t_n}(q_n,q)=A^{t_n-t}(q_n,\theta_n)+A^t(\theta_n,q)$.
This implies that 
\begin{align*}
\Phi^a(q_n,q)&\geq A^{t_n}(q_n,q)+at_n-1/n\\
&\geq A^{t_n-t}(q_n,\theta_n)+a(t_n-t)+A^t(\theta_n,q)
+at-1/n\\
&\geq  \Phi^a(q_n,\theta_n)+A^t(\theta_n,q)+at-1/n.
\end{align*}
This inequality implies that 
$$
u_n(q)\geq u_n(\theta_n)+A^t(\theta_n,q)+at -1/n.
$$
Since the functions $u_n$ are equi-Lipschitz, this implies
that the sequence $\theta_n$ is bounded,
by Lemma \ref{estimeA}.
By taking a subsequence, we assume that $\theta_n$ has a limit $\theta$, and, at the limit, we obtain
$$
u(q)\geq u(\theta)+A^t(\theta,q)+at,
$$
which implies that $u(q)\geq {\mathbf{T}}^tu(q)+ta$.
\qed

\section{Return to the periodic case.}
A more  precise link can be established between the contents of Sections \ref{WKAM} and \ref{regular} under the assumption that $H$ is periodic
(see Hypothesis \ref{h5}).
It is useful first to expose a variant of Section 
\ref{regular} adapted to the periodic case.
We leave  as exercises the proofs which are direct adaptations
of the ones given above.
From now on, we assume Hypotheses \ref{h1}, \ref{h2}, \ref{h3}, \ref{h4}, \ref{h5}.

\begin{thm}\label{C11p}
If (\ref{HJa}) admits a periodic sub-solution, then it admits a
periodic $C^{1,1}$ sub-solution.
Moreover, the set of periodic  $C^{1,1}$ sub-solutions is dense in the 
set of all periodic  sub-solutions for the uniform topology.
\end{thm}

\begin{defn}
The periodic critical value of $H$ is the real number $\alpha(0)$
defined as  the infimum of all real numbers 
$a$ such that (\ref{HJa}) has a periodic sub-solution.
The periodic  sub-solutions at level $\alpha(0)$ are called critical periodic sub-solutions. 
\end{defn}

\begin{defn}
The projected periodic Aubry set is the set $\mA(0)\subset \Tm^d$ of points 
$q$ such that the equality $H(q,du(q))=\alpha(0)$ holds for all
$C^1$ periodic  critical sub-solutions $u$.
\end{defn}

\begin{lem}
If $q\in \mA(0)$, then all $C^1$ critical periodic sub-solutions $u$
have the same differential at $q$. In other words, the restriction
$\Gamma_{u|\mA}$ does not depend on the $C^1$ critical 
periodic sub-solution 
$u$.
\end{lem}

\begin{prop}
There exists a $C^{1,1}$  periodic critical sub-solution $u_0$
such that $H(q,du_0(q))<\alpha(0)$
outside of $\mA(0)$.
\end{prop}

Without surprise,  we define the periodic Aubry set $\mA^*(0)$
as 
$$
\mA^*(0):= \Gamma_{u_0|\mA},
$$
with $u_0$ given by the proposition
(there is not a single $u_0$, but the Aubry set is well defined).

\begin{prop}
The set $\mA^*(0)\subset \Tm^d \times \Rm^{d*}$
is compact, non empty, and invariant.
\end{prop}

\proof
Let us prove that $\mA(0)$, hence $\mA^*(0)$ is not empty.
Assuming by contradiction that it was empty, then the equality
$H(q,du_0(q))<\alpha(0)$ would hold  for all $q\in \Rm^d$.
Since the function $q\lmto H(q,du_0(q))$ is periodic,
we could conclude that $\sup_q H(q,du_0(q))<\alpha(0)$,
which is in contradiction with the definition of $\alpha(0)$.
\qed

We are now in a position  to  specify  the connection with the invariant sets 
introduced in 
Section \ref{WKAM}:
\begin{prop}\label{egalite}
In the periodic case, we have the  equality
$$
\mA^*(0)=\cap_u \mI^*(u),
$$
where the intersection is taken on all periodic  weak KAM solutions.
\end{prop}

\proof
The inclusion $\mA^*(0)\subset \cap_u \mI^*(u)$
is proved as in Section \ref{regular}.
Our goal  is to prove the other inclusion.
Let $u_0$ be a $C^{1,1}$ periodic
sub-solution which is strict outside of $\mA(0)$.
The map $t\lmto {\mathbf{T}}^tu_0+t\alpha(0)$ is non-decreasing.
In addition, the functions ${\mathbf{T}}^tu_0+t\alpha(0)$ are equi-Lipschitz,
and they coincide with $u_0$ on $\mA$, hence they are equi-bounded.
As a consequence, 
${\mathbf{T}}^tu_0+t\alpha \lto u_{\infty}$ uniformly as $t\lto \infty$.

\textbf{Claim: } The limit  $u_{\infty}$ is a periodic  weak KAM solution such that 
$u_0< u_{\infty}$
 outside of $\mA(0)$.

In order to prove  that $u_{\infty}$ is a weak KAM solution, it is enough to notice that the function ${\mathbf{T}}^{t+s}u_0+(t+s)\alpha(0)$
converges both to $u_{\infty}$ and to ${\mathbf{T}}^su_{\infty}+s\alpha(0)$
when $t\lto \infty$. This implies, as desired, that 
${\mathbf{T}}^su_{\infty}+s\alpha(0)=u_{\infty}$ for each $s\geq 0$.

We know that $u_{\infty}\geq u_0$, with equality on $\mA(0)$.
 Conversely, let us consider a point $q$ 
such that $u_{\infty}(q)=u_0(q)$. The point $q$ is 
minimizing the difference 
$u_{\infty}-u_0$. Since $u_{\infty}$ is semi-concave and $u_0$
is $C^1$, the function $u_{\infty}$ must be differentiable at
$q$ with $du_{\infty}(q)=du_0(q)$. Since $u_{\infty}$
solves the Hamilton-Jacobi equation at its points
of differentiability, we conclude that 
$H(q,du_0(q))=H(q,du_{\infty}(q))=\alpha(0)
$,
hence $q\in \mA(0)$. We have proved the claim.

Let us now establish that $\mI(u_{\infty})=\mA(0)$,
which implies the proposition.
By Lemma $\ref{Tu=u}$, we have $\check {\mathbf{T}}^tu_{\infty}-t\alpha=u_{\infty}$ on $\mI(u_{\infty})$ for each $t\geq 0$.
Setting $\epsilon(t)=\sup(u_{\infty}-{\mathbf{T}}^tu_0-t\alpha(0))$,
we have 
$$
u_{\infty}\geq u_0
\geq \check {\mathbf{T}}^t \circ {\mathbf{T}}^tu_0\geq 
\check {\mathbf{T}}^t(u_{\infty}-\epsilon(t)-t\alpha(0) )
\geq \check {\mathbf{T}}^tu_{\infty}-\epsilon(t)-t\alpha(0)
= u_{\infty}-\epsilon(t)
$$
on $\mI(u_{\infty})$. Since this holds for all $t\geq 0$,
and since $\lim_{t\lto \infty} \epsilon(t)=0$, we conclude
that $u_0=u_{\infty}$ on $\mI(u_{\infty})$.
On the other hand, 
 we have seen that 
$
 u_0 <u_{\infty}
$
outside of $\mA(0)$, hence 
 $\mI(u_{\infty})\subset \mA(0)$.
\qed

We finish with an easy remark which is specific
to the periodic case:

\begin{prop}
All periodic weak KAM solutions have level $\alpha(0)$.
\end{prop}

\proof
Let $u_0$ be a critical periodic sub-solution, 
and let $u$ be a periodic  weak KAM solution at level $a$.
Since $u$ is a periodic sub-solution at level $a$,
the definition of $\alpha(0)$ implies that $a\geq \alpha(0)$.
On the other hand, there exists a constant $C$ such that 
$u-C\leq u_0\leq u+C$, which implies
$$
u={\mathbf{T}}^tu+ta\geq {\mathbf{T}}^tu_0-C+ta
\geq u_0+t(a-\alpha(0))-C\geq u+t(a-\alpha(0))-2C.
$$
We obtain  that $t(a-\alpha(0))\leq 2C$ for each $t\geq 0$, hence 
$a-\alpha(0)\leq 0$.
\qed

\section{The Lagrangian.}
In most expositions of  weak KAM theory, the Lagrangian plays an 
important role. In the present section, we 
relate it to  our main objects
in order to facilitate the connection with the core of the literature,
where  what we state here as properties is usually
taken as definitions.
We define the Lagrangian as
\begin{align*}
L:\Rm\times \Rm^d\times \Rm^d &\lto \Rm\\
(t,q,v) &\lmto \sup_{p\in (\Rm^d)^* }
\big(p\cdot v-H(t,q,p)\big).
\end{align*}
By standard results on convex Analysis, see\cite{Ro} for example,  we then have
$$
H(t,q,p)=\sup_{v\in \Rm^d}
 \big(p\cdot v-L(t,q,v)\big).
$$
We obviously have the Legendre inequality
$$
H(t,q,p)+L(t,q,v)\geq p\cdot v
$$
for all $t,q,p,v$.
This inequality is an equality if and only if 
$$
p=\partial_vL(t,q,v) \text{ or equivalently }
v=\partial_pH(t,q,p).
$$
Let $q(t):]t_0,t_1[\lto M$ be a curve,
The \textbf{action} of $q$ is the number
$$
\int_{t_0}^{t_1} L(t,q(t),\dot q(t))dt. 
$$
We can also call it Lagrangian action if we want to distinguish
from the previously defined Hamiltonian action.
The Lagrangian and Hamiltonian actions
are related as follows:

The Hamiltonian action of a curve $(q(t),p(t))$
is smaller than the Lagrangian action of its projection $q(t)$,
with equality if and only if 
$p(t)\equiv \partial_vL(t,q(t),\dot q(t))$.
In particular, the Hamiltonian action of an orbit is equal
to the Lagrangian action of its projection.

\begin{lem}
let $q_0$ and $q_1$ be two points of $\Rm^d$, and $t_0, t_1$
be two times, with $0<t_1-t_0<\sigma$. If $(q(s),p(s))$
is the orbit satisfying $q(t_0)=q_0, q(t_1)=q_1$, we have 
$$
S_{t_0}^{t_1}(q_0,q_1)
=\int_{t_0}^{t_1} L(s,q(s),\dot q(s))ds=\min_{\theta(s)} \int_{t_0}^{t_1}
 L(s,\theta(s),\dot \theta(s))ds,
$$
where the minimum is taken on the set of Lipschitz curves 
$\theta:[t_0,t_1]\lto \Rm^d$
which satisfy $\theta(t_0)=q_0$ and $\theta(t_1)=q_1$.
\end{lem}

\proof
Since $S_{t_0}^{t_1}(q_0,q_1)$ is the Hamiltonian 
action of the unique orbit $(q(t),p(t))$, it is also the Lagrangian
action of the curve $q(t)$:
$$
S_{t_0}^{t_1}(q_0,q_1)
= \int_{t_0}^{t_1} L(s,q(s),\dot q(s))ds.
$$
The function $u(t,q):= S_{t_0}^t(q_0,q)$
solves (\ref{HJ}) on $]t_0,t_1[\times \Rm^d$.
Let us now consider any Lipschitz curve 
$\theta(s):[t_0,t_1]\lto \Rm^d$
satisfying $\theta(t_0)=q_0$ and $\theta(t_1)=q_1$,
and write
\begin{align*}
\makebox[2.5cm]{}
\int_{t_0}^{t_1}
 L(s,\theta(s),\dot \theta(s))ds
& \geq \int_{t_0}^{t_1} \partial_qu(s,\theta(s))\cdot \dot \theta(s)-
H(s,\theta(s),\partial_qu(s,\theta(s)) ds\\
& =\int_{t_0}^{t_1} \partial_qu(s,\theta(s))
 \cdot \dot \theta(s)-
\partial_t u(s,\theta(s))ds\\
&= u(t_1,q_1)-u(t_0,q_0)=S_{t_0}^{t_1}(q_0,q_1).
\makebox[4cm][r]{\fbox{}}
\end{align*}

The following proposition is usually taken as the definition of $A$:
\begin{prop}
Given two points  $q_0$ and $q_1$  and two times  $t_0< t_1$, we have
$$
A_{t_0}^{t_1}(q_0,q_1)=\min_{\theta(s)}\int_{t_0}^{t_1}
 L(s,\theta(s),\dot \theta(s))ds,
$$
where the minimum is taken on the set of Lipschitz curves 
$\theta:[t_0,t_1]\lto \Rm^d$
which satisfy $\theta(t_0)=q_0$ and $\theta(t_1)=q_1$.
\end{prop}
It is part of the statement that the minimum is achieved. This 
is usually called the Theorem of Tonelli.
The statement can be extended to absolutely continuous curves
instead of Lipschitz curves, but this setting is not useful 
for our discussion.

\proof
For $n$ large enough, we have $(t_1-t_0)/n<\sigma$,
hence, setting $\tau_i=t_0+i(t_1-t_0)/n$, 
\begin{align*}
\makebox[1cm]{}
A_{t_0}^{t_1}(q_0,q_1)&=\min_{(\theta_1,\ldots,\theta_{n-1})}\big(
S_{t_0}^{\tau_1}(q_0,\theta_1)+S_{\tau_1}^{\tau_2}(\theta_1,\theta_2)+
\cdots+
S_{\tau_{n-1}}^{t_1}(\theta_{n-1},q_1)
\big)\\
&=
\min_{(\theta_1,\ldots,\theta_{n-1})}\big(
\min_{\theta(s)}\int_{t_0}^{\tau_1}
L(s,\theta(s),\dot \theta(s))
ds
+\ldots +
\min_{\theta(s)}\int_{\tau_{n-1}}^{t_1}
L(s,\theta(s),\dot \theta(s))
ds
\big)\\
&=\min_{\theta(s)} \int_{t_0}^{t_1}L(s,\theta(s),\dot \theta(s))
ds.
\makebox[8cm][r]{\fbox{}}
\end{align*}

\begin{appendix}

\section{Some technical results.}

\begin{prop}\label{global}
A Lipschitz   map $F:\Rm^d\lto \Rm^d$ which satisfies
$Lip(F-Id)<1$ is a bi-Lipschitz homeomorphism of $\Rm^d$.
Its inverse is Lipschitz and $Lip(F^{-1})\leq (1-k)^{-1}$. If $F$ is $C^1$, then so is $F^{-1}$.
\end{prop}

\proof
The equation 
$F(q)=\theta$ can be rewritten
$$
\theta-(F(q)-q)=q
$$
The map on the left being contracting, we conclude that $F$ is invertible.
We now write 
$$
|x_1-x_0|-|F(x_1)-F(x_0)|\leq |(F(x_1)-x_1)-(F(x_0)-x_0)|\leq k|x_1-x_0|
$$
and deduce that $|F(x_1)-F(x_2)|\geq (1-k) |x_1-x_0|$.
\qed

\begin{prop}\label{monotonemap}
Let $F:\Rm^d\lto \Rm^d$ be a $C^1$, $c$-monotone map 
on $\Rm^d$, with $c>0$. Then $F$ is a diffeomorphism
from $\Rm^d$ onto itself.
\end{prop}

\proof
Let us consider a point $\theta\in\Rm^d$, and the 
line 
$\theta(s)=F(0)+s(\theta-F(0))$.
Since $F$ is a local diffeomorphism around $0$,
the points $\theta(s)$ for small $s$ have a unique preimage $p(s)$.
Let $S$ be the infimum of the positive real numbers $s$ such that 
the equation $F(p)=\theta(s)$ does not have a solution in 
$\Rm^d$.
 The curve $p(s)$ is well-defined, $C^1$, and Lipschitz
on $[0,S[$, hence, if $S$ is finite, it extends at $S$ with 
with  $F(p(S))=\theta(S)$.
Since $F$ is a local diffeomorphism at $p(S)$,
the points near $\theta(S)$ have preimages, which
contradicts the definition of $S$.
Hence $S$ can't be finite.
\qed

\begin{lem}\label{a/b}
Let $A$ be a $d\times d$ matrix, such that $A\geq aId$
in the sense of quadratic forms, and $\|A\|\leq b$.
Then $A^{-1}\geq (a/b^2)I$ in the sense of quadratic forms.
\end{lem}

\proof
We have 
$$
(A^{-1}v,v)=(A A^{-1}v,A^{-1}v)\geq a|A^{-1}v|^2\geq
a(|v|/b)^2.
$$
\qed

The following important result appears in Fathi's book on Weak KAM theory
(the proof is also his):

\begin{prop}\label{fathi}
Let $u:\Rm^d\lto \Rm$ be a  function and $K$
be a positive number.
Let $\mI\in \Rm^d$
be the set of points where  $u$ has both a $K$-super-differential
and a $K$-sub-differential.
Then, the function $u$ is differentiable at each point of 
$\mI$ and the function $q\lmto du(q)$ is $6K$-Lipschitz on $\mI$.
\end{prop}

\proof
For each $q\in \mI$, there exists a unique  $l(q)\in \Rm^{d*}$ such that 
$$|u(q+\theta)-u(q)-l(q)\cdot \theta|\leq K\|\theta\|^2.$$
We conclude that $l(q)$ is the differential of $u$ at $q$, and we have to prove that the map
$q\lmto l(q)$ is Lipschitz on $\mI$.
We have, for each $q$, $\theta$ and $y$ in $H$:
$$
l(q)\cdot (y+\theta)-K\|y+\theta\|^2\leq u(q+y+\theta)-u(q)\leq l(q)\cdot (y+\theta)+K\|y+\theta\|^2
$$
$$
l(q+y)\cdot (-y)-K\|y\|^2\leq u(q)-u(q+y)\leq
l(q+y)\cdot (-y)+K\|y\|^2
$$
$$
l(q+y)\cdot (-\theta)- K\|\theta\|^2 \leq u(q+y)-u(q+y+\theta)\leq l(q+y)\cdot (-\theta) +K\|\theta\|^2.
$$
Taking the sum, we obtain
$$
\big|(l(q+y)-l(q))\cdot (y+\theta)
\big|\leq K\|y+\theta\|^2+K\|y\|^2+K\|z\|^2.
$$
By a change of variables, we get 
$$
\big|(l(q+y)-l(q))\cdot \theta
\big|\leq K\|\theta\|^2+K\|y\|^2+K\|\theta-y\|^2.
$$
Taking $\|\theta\|=\|y\|$, we obtain
$$
\big|(l(q+y)-l_(q))\cdot (\theta)
\big|\leq 6K\|\theta\|\|y\|
$$
for each $\theta$ such that $\|\theta\|=\|y\|$,
we conclude that 
$$\|l(q+y)-l(q)\|\leq 6K\|y\|.
$$
\qed

\begin{lem}\label{infimum}
Let $u$ be a finite valued function  which is the infimum
of a family $\mF$ of equi-semi-concave functions:
$u=\inf_{f\in \mF} f$. Then the function $u$ is 
semi-concave.
\end{lem}
It is important in the statement to assume that $u$
is really finite valued at each point.

\proof
Let us assume that the functions in $\mF$
are $k$-semi-concave.
Given a point $q_0\in \Rm^d$ let 
$f_n(q)=a_n+p_n\cdot q +k/2\|q\|^2$
be a sequence of functions of $\mF$
such that $f_n(q_0)\lto u(q_0)$.
We have 
$f_n(q)\leq f_n(q_0)+p_n\cdot (q-q_0)+k/2\|q-q_0\|^2$
for some sequence $p_n\in \Rm^*$.
If the sequence $p_n$ is bounded, 
then we can take the limit along a subsequence 
and get the inequality
$$
u(q)\leq u(q_0)+p\cdot (q-q_0)+k/2\|q-q_0\|^2.
$$
If this holds for each $q_0$, we conclude that 
$u$ is $k$-semi-concave.
Let us now prove that $p_n$ is bounded.
If this was not true, there would exist a point $q$
such that $p_n\cdot (q-q_0)$ is not bounded 
from below. This would imply that 
$$
u(q)=\inf _{f\in \mF}f(q)\leq 
\inf_n f_n(q)\leq \inf_n 
\big(f_n(q_0)+p_n\cdot (q-q_0)+k/2\|q-q_0\|^2
\big)
=
-\infty,
$$
which would contradict the finiteness of $u$ at $q$.
\qed

\textbf{Acknowledgements. }
Many  thanks  to Lyonell Boulton and Sergei Kuksin for organizing this CANPDE session on Weak KAM theory.
I also thank Albert Fathi form many useful comments and suggestions.

\end{appendix}
\begin{small}

\end{small}
\end{document}